\documentclass[11pt,twoside]{article}
\usepackage[hmargin=0.8in,vmargin=0.9in]{geometry}
\usepackage{fancyhdr}
\usepackage{graphicx}
\usepackage{subfigure}
\usepackage{amssymb}
\usepackage{amsmath}
\usepackage{amsfonts}
\usepackage{adjustbox}
\usepackage{theorem}
\usepackage{mathrsfs}
\usepackage{mathtools}
\usepackage{bm}
\usepackage{tabularx}
\usepackage{color}
\usepackage{setspace}
\usepackage{exscale}
\usepackage{relsize}
\usepackage{epstopdf}
\usepackage{float}
\usepackage{cite}

\usepackage{hyperref}
\usepackage{setspace}

\usepackage{algorithm}
\usepackage{algpseudocode}

\usepackage{color}
\usepackage{dsfont}

\usepackage{booktabs,multirow} % for much better looking tables
\usepackage{array} % for better arrays (eeg matrices) in maths
\usepackage{paralist} % very flexible & customisable lists (eeg. enumerate/itemize, etc.)
\usepackage{verbatim} % adds environment for commenting out blocks of text & for better verbatim
\usepackage{subfigure} % make it possible to include more than one captioned figure/table in a single float
\theoremstyle{plain}			% use "default" font

{\theorembodyfont{\rmfamily}}

\newenvironment{DA}{{\flushleft \bf Declarations:}}{}
\newtheorem{rmk}{Remark}[section]

\allowdisplaybreaks[1]

\numberwithin{equation}{section}
\numberwithin{figure}{section}
\numberwithin{table}{section}

\newcommand\eref[1]{(\ref{#1})}

\newcommand*\xbar[1]{%
  \hbox{%
    \vbox{%
      \hrule height 0.5pt % The actual bar
      \kern0.4ex%         % Distance between bar and symbol
      \hbox{%
        \kern-0.05em%      % Shortening on the left side
        \ensuremath{#1}%
        \kern-0.05em%      % Shortening on the right side
      }%
    }%
  }%
}

\newcommand{\mF}{\bm{F}}

\newcommand{\mG}{\bm{G}}

\newcommand{\mU}{\bm{U}}
\newcommand{\bmu}{\bm{u}}

\newcommand{\dx}{\Delta x}
\newcommand{\dy}{\Delta y}

\newcommand{\hf}{{\frac{1}{2}}}

\newcommand{\jph}{{j+\frac{1}{2}}}
\newcommand{\jmh}{{j-\frac{1}{2}}}
\newcommand{\iph}{{i+\frac{1}{2}}}
\newcommand{\imh}{{i-\frac{1}{2}}}
\newcommand{\kph}{{k+\frac{1}{2}}}
\newcommand{\kmh}{{k-\frac{1}{2}}}

\def\softd{{\leavevmode\setbox1=\hbox{d}%
		\hbox to 1.05\wd1{d\kern-0.4ex{\char039}\hss}}}

\title{Numerical Realization of an Entropy-Based Selection Principle in the Example of Kelvin-Helmholtz Instability}
\author{Shaoshuai Chu\thanks{Department of Mathematics, RWTH Aachen University, 52056 Aachen, Germany; {\tt chu@igpm.rwth-aachen.de}}~ and
~Michael Herty\thanks{Department of Mathematics, RWTH Aachen University, 52056 Aachen, Germany; Department of Mathematics and Applied
Mathematics, University of Pretoria, Hatfield, 0028, South Africa; {\tt herty@igpm.rwth-aachen.de}}}
\date{}

\begin{document}
\maketitle
\begin{abstract}
In this paper, we use Young measures to incorporate an entropy-based selection criterion into numerical fluxes for multidimensional compressible Euler flows. At each spatial degree of freedom, we solve a linear-programming (LP) problem to select a discrete probability measure whose first moment equals the evolved conservative state. The LP minimizes the expected mathematical entropy, equivalently maximizing the expected physical entropy among the feasible measures. The nonlinear Euler fluxes are then averaged with respect to a selected minimizer and used at each stage of the time discretization, thereby incorporating the prescribed criterion into the numerical evolution. To this end, we construct two-dimensional Young-measure schemes with first-, second-, third-, fifth-, seventh-, and ninth-order spatial reconstructions and apply them to the Kelvin--Helmholtz instability. We compare them with corresponding standard local Lax--Friedrichs (LLF) flux-splitting schemes and LLF schemes equipped with an entropy-based local characteristic decomposition (LLF--ELCD). The numerical results reveal scheme-dependent flow patterns, particularly in the small-scale vortical structures generated during the roll-up of the shear layers. Averaging over time and reconstruction orders substantially reduces the differences in the density profiles, while regional empirical density distributions retain scheme-dependent differences. At every reconstruction order considered, the Young-measure scheme yields the largest time-averaged domain-integrated physical entropy among the schemes of the same order. For the averages over all six orders, the Young-measure schemes also yield the smallest time-averaged energy defect and Bregman functional among the three scheme families. By contrast, the LLF--ELCD schemes do not systematically yield larger physical entropy than standard LLF. These results provide numerical evidence for an entropy preference consistent with the local LP objective and support the use of Young-measure closures to incorporate prescribed selection criteria into multidimensional flow computations.
\end{abstract}



\section{Introduction}
In this paper, we consider the two-dimensional (2-D) Euler equations of gas dynamics:
\begin{equation}
\begin{aligned}
&\rho_t+(\rho u)_x+(\rho v)_y=0,\\
&(\rho u)_t+(\rho u^2+p)_x+(\rho uv)_y=0,\\
&(\rho v)_t+(\rho uv)_x+(\rho v^2+p)_y=0,\\
&E_t+\left[u(E+p)\right]_x+\left[v(E+p)\right]_y=0.
\end{aligned}
\label{1.1}
\end{equation}
Here, $x$ and $y$ are the spatial variables, $t$ is time, $\rho$ is the density, $u$ and $v$ are the velocities in the $x$- and $y$-directions, respectively, $E$ is the total
energy density, and $p$ is the pressure. The system \eref{1.1} is closed by the polytropic equation of state
\begin{equation}
e(\rho,\theta)=c_v\theta,\qquad
c_v=\frac{1}{\gamma-1},\qquad
p(\rho,\theta)=\rho\theta=(\gamma-1)\rho e
=(\gamma-1)\left(E-\frac{|\bm m|^2}{2\rho}\right),
\label{1.2}
\end{equation}
where $e$ is the specific internal energy, $\theta$ is the absolute temperature, $\bm m=(\rho u,\rho v)^\top$ is the momentum vector, and $\gamma>1$ is the ratio of specific
heats. The second law of thermodynamics is expressed through
\begin{equation}
S_t+(Su)_x+(Sv)_y\geq0,\qquad
S:=c_v\rho\ln\left(\frac{p}{\rho^\gamma}\right),
\label{1.3}
\end{equation}
where $S$ is the physical entropy density. Correspondingly, $\eta:=-S$ is the mathematical entropy, for which the inequality in \eref{1.3} is reversed. On a bounded domain $\Omega\subset\mathbb R^2$ with periodic or no-flux boundary conditions, we prescribe
\begin{equation}
\rho(x,y,0)=\rho_0(x,y),\qquad
\bm m(x,y,0)=\bm m_0(x,y),\qquad
E(x,y,0)=E_0(x,y).
\label{1.4}
\end{equation}
The corresponding initial pressure and physical entropy are
$$
p_0=(\gamma-1)\left(E_0-\frac{|\bm m_0|^2}{2\rho_0}\right),
\qquad
S_0=c_v\rho_0\ln\left(\frac{p_0}{\rho_0^\gamma}\right).
$$

Solutions of \eref{1.1}--\eref{1.2} may develop shocks, rarefaction waves, and contact discontinuities even from smooth initial data and must then be understood in the weak sense. In one space dimension, the entropy inequality provides an important admissibility criterion. In several space dimensions, however, entropy admissibility alone does not guarantee uniqueness. Convex-integration results show that multidimensional compressible Euler systems may admit infinitely many admissible weak solutions for certain initial data; see, e.g., \cite{de2010admissibility,chiodaroli2015global,feireisl2020density} and the references therein.

These nonuniqueness results motivate the study of generalized solutions and additional selection criteria. For oscillatory numerical approximations, the weak limit of the conservative variables does not, in general, determine the limits of the nonlinear Euler fluxes. Young measures retain the local distribution of states and therefore provide a framework for representing such oscillations; see \cite{diperna1985measure,fjordholm2016computation}. For compressible flows, dissipative weak (DW) solutions account for oscillation and concentration effects through Young measures and compatible energy and Reynolds stress defects. Their relation to consistent numerical approximations is developed in \cite{Lukacova_book}. Averaging provides a way to study the limiting behavior of these approximations. Under suitable consistency and boundedness assumptions, Ces\`aro averages of appropriate subsequences yield strong convergence of the density, momentum, and entropy to the corresponding variables of a DW solution; see, e.g., \cite{Feireisl2020_Kconvergence,feireisl2021computing,
Lukacova_book,lukavcova2023limit}. Related numerical studies have used averaging and statistical diagnostics to investigate deterministic and random Euler flows. In \cite{CHIIKL}, stochastic collocation was combined with a fifth-order alternative weighted essentially non-oscillatory (A-WENO) solver, and Ces\`aro averages over embedded spatial meshes were used to evaluate statistical quantities. The study in \cite{CHKLY26} compared several numerical methods using cumulative averages, empirical distributions, entropy, energy defects, and a Bregman functional. In that study, the selection criteria were used to assess already computed numerical solutions.

In the present work, we study a Young-measure closure in which a prescribed selection criterion enters the flux evaluation. To this end, we extend the LP-based Young-measure schemes developed in \cite{CH_high_order,CHLZ25} from one to two space dimensions. At each spatial degree of freedom, we construct a discrete probability measure whose first moment equals the evolved conservative state. The normalization and first-moment constraints may admit several measures, which can yield different averages of the nonlinear Euler fluxes. The objective function of the linear-programming (LP) problem is used to select a measure among those satisfying these constraints. We use the mathematical entropy $\eta=-S$ as the LP objective. For each prescribed conservative state, minimizing its expectation is equivalent to maximizing the expected physical entropy among the feasible discrete measures. The Euler fluxes are then averaged with respect to a selected minimizer and used in the numerical update. In this way, the entropy-based criterion enters the numerical evolution through the fluxes evaluated at every stage of the time discretization. We apply the resulting 2-D schemes to the Kelvin--Helmholtz (KH) instability. As the shear layers roll up, increasingly fine vortical structures may appear, and the computed flow patterns can depend strongly on the numerical method and mesh resolution; see, e.g., \cite{CCHKL_22,CCK23_Adaptive,feireisl2021computing,CHKLY26}. This sensitivity provides a useful setting for examining how the prescribed selection criterion influences the computed evolution.

To determine whether this explicit optimization produces an observable entropy preference, we compare three families of schemes: the proposed Young-measure schemes with entropy-based LP closure, standard local Lax--Friedrichs (LLF) flux-splitting schemes, and LLF schemes incorporating the entropy-based local characteristic decomposition (ELCD) from \cite{CHK_NEW_LCD}. For each family, we construct first-, second-, third-, fifth-, seventh-, and ninth-order realizations. This allows comparisons between methods of the same formal order. Since the first-order method requires no LCD, the standard LLF and LLF--ELCD schemes coincide at this order. The third family is included, in particular, to determine whether incorporating entropy into a classical LCD is sufficient to produce an entropy-selective numerical evolution.  The numerical experiments show that the three families of schemes may produce different evolutions of the KH instability, with the clearest differences appearing in the small-scale vortical structures. The density profiles produced by standard LLF and LLF--ELCD are quite similar, whereas noticeable differences arise between these profiles and the Young-measure results, especially at higher orders. We also examine finite cumulative averages over different reconstruction orders on a fixed mesh. The later averages show closer agreement, and the differences between the three families are substantially smaller after time averaging. The density profiles averaged over both time and reconstruction orders are remarkably similar for the three families, while their regional empirical density distributions retain scheme-dependent differences. 

At every formal order considered, the Young-measure scheme yields the largest time-averaged domain-integrated physical entropy among the corresponding schemes. This systematic preference is consistent with the imposed objective and provides numerical evidence that the LP closure affects the selected evolution. By contrast, the LLF--ELCD schemes do not systematically yield larger entropy than standard LLF. Thus, using entropy in the LCD alone does not ensure the entropy preference observed for the Young-measure schemes. This comparison indicates that the explicit optimization in the LP problem plays an important role in the resulting numerical evolution. Moreover, for the averages over all six orders, the Young-measure schemes yield the largest time-averaged total physical entropy, the smallest time-averaged energy defect, and the smallest time-averaged Bregman functional among the three scheme families.  By contrast, the LLF--ELCD schemes yield the largest time-averaged energy defect. Therefore, the numerical approximation selected from the three scheme families depends on the prescribed selection criterion. 

These findings indicate that the objective function is an active modeling and selection component rather than an auxiliary closure tool. This interpretation points to a broader use of the framework. Although we focus on mathematical entropy, the formulation is not intrinsically restricted to this choice: other objective functions may be introduced to encode alternative physically or mathematically motivated criteria while maintaining the admissibility constraints. The presented framework thus provides a flexible and systematic mechanism for numerical approximations consistent with a prescribed criterion. Optimizing a selected quantity, however, does not by itself imply greater numerical accuracy. In the present experiments, the larger physical entropy shows that the Young-measure approximations are more strongly aligned with the maximum-entropy criterion imposed through the LP.

The rest of the paper is organized as follows. In \S\ref{sec3}, we review the 2-D Young-measure, LLF, and LLF--ELCD schemes. In \S\ref{sec4}, we present the numerical results for the KH instability and examine the resulting flow structures, entropy evolution, time and cross-order averages, and regional density distributions. We also compare the numerical solutions using several selection criteria and study the density marginals of the averaged Young measures. Finally, we give some concluding remarks in \S\ref{sec5}.

\section{Numerical Schemes}\label{sec3}
In this section, we present a unified framework for the first-, second-, third-, fifth-, seventh-, and ninth-order numerical schemes used in the numerical experiments. In particular, we describe the Young-measure schemes and the two families of classical LLF flux-splitting schemes used for comparison.

We first rewrite the 2-D Euler equations of gas dynamics \eref{1.1} in the vector form
\begin{equation*}
\mU_t+\mF(\mU)_x+\mG(\mU)_y=\bm 0,
\end{equation*}
where $\mU=(\rho,\rho u,\rho v,E)^\top$, $\mF(\mU)=\big(\rho u,\rho u^2+p,\rho uv,u(E+p)\big)^\top$, and $\mG(\mU)=\big(\rho v,\rho uv,\rho v^2+p,v(E+p)\big)^\top$.
Assume that $(x_j,y_k)$ are the points of a uniform Cartesian grid satisfying $x_{j+1}-x_j\equiv\dx$ and $y_{k+1}-y_k\equiv\dy$. We denote the corresponding numerical solution by $\mU_{j,k}(t)\approx\mU(x_j,y_k,t)$. The point values $\mU_{j,k}$ are evolved in time by numerically solving the following system of ODEs:
\begin{equation}\label{equ3.2}
\frac{{\rm d}\mU_{j,k}}{{\rm d}t}=-\frac{\bm{{\cal F}}^{(r)}_{\jph,k}-\bm{{\cal F}}^{(r)}_{\jmh,k}}{\dx}-\frac{\bm{{\cal G}}^{(r)}_{j,\kph}-\bm{{\cal G}}^{(r)}_{j,\kmh}}{\dy},
\end{equation}
where $r\in\{1,2,3,5,7,9\}$ is the prescribed order of spatial accuracy. Note that all of the indexed quantities are time-dependent, but from here on, we will suppress this dependence for the sake of brevity.

The numerical fluxes $\bm{{\cal F}}^{(r)}_{\jph,k}$ and $\bm{{\cal G}}^{(r)}_{j,\kph}$ are constructed dimension by dimension \cite{Jiang96}. More precisely, the flux values are first split using the LLF flux splitting and are then reconstructed at the corresponding cell interfaces. Let ${\cal R}^{(r),+}$ and ${\cal R}^{(r),-}$ denote the left- and right-biased reconstruction operators of order $r$, respectively. The numerical fluxes are defined by
\begin{equation}\label{equ3.3}
\begin{aligned}
\bm{{\cal F}}^{(r)}_{\jph,k}&={\cal R}^{(r),+}_{\jph,k}\bigg(\left\{\bm F_{\ell,k}^{+}\right\}_{\ell\in{\cal S}^{(r)}_{\jph}}\bigg)+{\cal R}^{(r),-}_{\jph,k}
\bigg(\left\{\bm F_{\ell,k}^{-}\right\}_{\ell\in{\cal S}^{(r)}_{\jph}}\bigg),\\
\bm{{\cal G}}^{(r)}_{j,\kph}&={\cal R}^{(r),+}_{j,\kph}\bigg(\left\{\bm G_{j,\ell}^{+}\right\}_{\ell\in{\cal S}^{(r)}_{\kph}}\bigg)+{\cal R}^{(r),-}_{j,\kph}
\bigg(\left\{\bm G_{j,\ell}^{-}\right\}_{\ell\in{\cal S}^{(r)}_{\kph}}\bigg).
\end{aligned}
\end{equation}
Here, the split fluxes are
\begin{equation}\label{equ3.4}
\bm F_{\ell,k}^{\pm}=\frac{1}{2}\bigg(\bm\Phi_{\ell,k}\pm a^x_{\jph,k}\mU_{\ell,k}\bigg),\quad \ell\in{\cal S}^{(r)}_{\jph} \quad {\rm and} \quad 
\bm G_{j,\ell}^{\pm}=\frac{1}{2}\bigg(\bm\Psi_{j,\ell}\pm a^y_{j,\kph}\mU_{j,\ell}\bigg),\quad \ell\in{\cal S}^{(r)}_{\kph},
\end{equation}
where $\bm\Phi_{j,k}$ and $\bm\Psi_{j,k}$ are the unsplit flux values, which will be specified below. The local speeds in \eref{equ3.4} are estimated by
\begin{equation*}
a^x_{\jph,k}=\max_{\ell\in{\cal S}^{(r)}_{\jph}}\big(|u_{\ell,k}|+c_{\ell,k}\big),\qquad a^y_{j,\kph}=\max_{\ell\in{\cal S}^{(r)}_{\kph}}\big(|v_{j,\ell}|+c_{j,\ell}\big),
\end{equation*}
where the speed of sound is determined by $c_{j,k}=\sqrt{\frac{\gamma p_{j,k}}{\rho_{j,k}}}$. The corresponding reconstruction stencils are
\begin{equation*}
{\cal S}^{(r)}_{\jph}:=\begin{cases}\{j,j+1\},& r=1,\\[1mm]
\{j-1,j,j+1,j+2\},& r=2,\\[1mm]
\{j-s+1,\ldots,j+s\},& r=2s-1,\quad s=2,3,4,5,
\end{cases}
\end{equation*}
and ${\cal S}^{(r)}_{\kph}$ is defined analogously in the $y$-direction.

For systems of conservation laws, the reconstructions of order $r\geq2$ are performed in local characteristic variables, a strategy that is often used in high-order schemes to reduce spurious oscillations and can also enhance the resolution of second-order schemes; see, e.g., \cite{Jiang13,Liu17,Nonomura20,Qiu02,Shu20} and references therein. In particular, each operator in \eref{equ3.3} includes the projection of all split fluxes on the stencil onto a single characteristic basis constructed at the corresponding interface, the componentwise scalar reconstruction in the characteristic variables, and the transformation of the reconstructed fluxes back to the conservative variables. The detailed formulas are provided in Appendix~\ref{appa}.

\subsection{Young-Measure Approach}
In the schemes using Young measures, a discrete measure is reconstructed at every spatial grid point and at every evaluation of the right-hand side of \eref{equ3.2}. Let $\bm z=(z^{(1)},\ldots,z^{(d)})^\top\in\mathbb R^d$, with $d=4$ in the 2-D Euler equations of gas dynamics, denote a generic phase-space variable, namely, a possible value of the conservative vector $\mU$. The phase-space grid points are denoted by $\bm z_\ell\in\mathbb R^d$, $\ell=1,\ldots,N_{\bmu}$. Here, $N_{\bmu}$ denotes the total number of phase-space grid points. If $N_q$ grid points are used in the $q$-th phase-space direction before the admissibility restriction described below, then the underlying tensor-product grid contains $\prod_{q=1}^dN_q$ points. If a uniform tensor-product grid is used in phase space, its cell volume is $\Delta V=\prod_{q=1}^d\Delta z^{(q)}$. We assume that the phase-space discretization is chosen so that the constraints introduced below are feasible; in particular, every evolved state $\mU_{j,k}$ must belong to the convex hull of the phase-space grid points.

At each grid point $(x_j,y_k)$, the values $\{\mu_{j,k,\ell}^*\}_{\ell=1}^{N_{\bmu}}$ are obtained by solving the LP problem
\begin{align}
\{\mu_{j,k,\ell}^*\}_{\ell=1}^{N_{\bmu}}&\in\operatorname*{argmin}_{\{\mu_\ell\}_{\ell=1}^{N_{\bmu}}}\Delta V\sum_{\ell=1}^{N_{\bmu}}\eta(\bm z_\ell)\mu_\ell,\label{linprog-a}\\
\mbox{subject to}\qquad &\mu_\ell\geq0,\qquad \ell=1,\ldots,N_{\bmu},\label{linprog-b}\\
&\mu_\ell\leq\frac{1}{\Delta V},\qquad \ell=1,\ldots,N_{\bmu},\label{linprog-c}\\
&\Delta V\sum_{\ell=1}^{N_{\bmu}}\mu_\ell=1,\label{linprog-d}\\
&\Delta V\sum_{\ell=1}^{N_{\bmu}}\bm z_\ell\mu_\ell=\mU_{j,k}.\label{linprog-e}
\end{align}
Thus, $\nu_{j,k}^*=\Delta V\sum_{\ell=1}^{N_{\bmu}}\mu_{j,k,\ell}^*\,\delta_{\bm z_\ell}$ is a discrete probability measure whose first moment $\Delta V\sum_{\ell=1}^{N_{\bmu}}\bm z_\ell\mu^*_{j,k,\ell}$ coincides with the evolved conservative state $\mU_{j,k}$. Notice that the upper bound in \eref{linprog-c} follows from the nonnegativity and normalization constraints \eref{linprog-b} and \eref{linprog-d}; it is retained here to match the LP problem used in the implementation.

In this paper, the objective function in \eref{linprog-a} is the mathematical entropy, and consequently, for the prescribed first moment $\mU_{j,k}$, minimizing the objective function in \eref{linprog-a} is equivalent to maximizing the expected physical entropy over all admissible discrete measures satisfying the constraints \eref{linprog-b}--\eref{linprog-e}. After the minimizer $\mu^*_{j,k,\ell}$ has been determined, the flux values in \eref{equ3.4} are computed by averaging the physical fluxes with respect to $\mu_{j,k,\ell}^*$:
\begin{equation*}
\bm\Phi_{j,k}^{\rm YM}:=\left\langle\nu_{j,k}^*,\mF\right\rangle=\Delta V\sum_{\ell=1}^{N_{\bmu}}\mF(\bm z_\ell)\mu_{j,k,\ell}^*,\quad
\bm\Psi_{j,k}^{\rm YM}:=\left\langle\nu_{j,k}^*,\mG\right\rangle=\Delta V\sum_{\ell=1}^{N_{\bmu}}\mG(\bm z_\ell)\mu_{j,k,\ell}^*.
\end{equation*}

For the Young-measure schemes, we set $\bm\Phi_{j,k}=\bm\Phi_{j,k}^{\rm YM}$ and $\bm\Psi_{j,k}=\bm\Psi_{j,k}^{\rm YM}$ in \eref{equ3.4}. For $r\geq2$, the characteristic basis used in the Young-measure schemes is constructed from the arithmetic averages of the primitive variables, as described in Appendix~\ref{appa}. The left and right eigenvector matrices used for the LCD are therefore not averaged with respect to the Young measure.

\subsection{Classical LLF Flux-Splitting Schemes}
For the classical LLF flux-splitting schemes, the unsplit flux values in \eref{equ3.4} are evaluated directly by the evolved conservative state:
\begin{equation*}
\bm\Phi_{j,k}=\mF(\mU_{j,k}),\qquad
\bm\Psi_{j,k}=\mG(\mU_{j,k}).
\end{equation*}
All other components of the discretization remain unchanged. We consider two choices of the characteristic basis. In the standard LLF schemes, the basis is constructed using the arithmetic averages of the primitive variables, as described in Appendix~\ref{appa}. In the LLF-ELCD schemes, the interface state used to construct the characteristic basis is selected from four candidate states according to their physical entropy; see Appendix~\ref{appb}.

It is important to distinguish the entropy-based LCD from the LP closure in \eref{linprog-a}--\eref{linprog-e}. The former uses entropy only to select an interface state for the LCD and does not optimize the numerical solution over a class of admissible measures. It therefore does not enforce a maximum-entropy selection principle. Moreover, since no LCD is needed when $r=1$, the two classical first-order LLF schemes coincide.

\begin{rmk}
The semi-discrete ODE system \eref{equ3.2} is integrated in time using the three-stage third-order strong stability preserving Runge-Kutta method (see, e.g., \cite{Gottlieb11,Gottlieb12}) for the second-, third-, fifth-, seventh-, and ninth-order schemes and the forward Euler method for the first-order scheme. For the Young-measure schemes, the LP closure is recomputed at every Runge-Kutta stage. In all of the numerical experiments, we take the CFL number to be $0.475$ and $\gamma=1.4$.  In the rest of the paper, we refer to the schemes of formal order $r\in\{1,2,3,5,7,9\}$ in the three families—the standard LLF flux-splitting schemes, the LLF flux-splitting schemes equipped with the entropy-based LCD, and the Young-measure schemes---as the $r$-Order LLF, $r$-Order LLF-ELCD, and $r$-Order YM schemes, respectively, for the sake of brevity.
\end{rmk}

\begin{rmk}
The LP problem \eref{linprog-a}--\eref{linprog-e} can be solved using the black-box solver \texttt{linprog} in MATLAB or the KKT-based solver described, for example,  in \cite[Appendix~D]{CH_high_order}. 
\end{rmk}

\section{Numerical Results}\label{sec4}
In this section, we apply the schemes described in \S\ref{sec3} to the KH instability problem. All computations are performed on the domain $[0,1]\times[0,1]$ using a uniform mesh with $\dx=\dy=1/512$.

To facilitate comparisons across formal order of the schemes, we use the following notation for the solutions computed within each scheme family, but with different orders of the resulting scheme: $(\rho_\ell,\bm m_\ell,S_\ell)$, $\ell=1,\dots,6$, are solutions computed by the first-, second-, third-, fifth-, seventh-, and ninth-order schemes, respectively, that is, $\ell(1)=1$, $\ell(2)=2$, $\ell(3)=3$, $\ell(5)=4$, $\ell(7)=5$, and $\ell(9)=6$. Motivated by the concept of $\mathcal{K}$-convergence described in
\cite{Feireisl2020_Kconvergence,feireisl2021computing}, we will average these solutions to obtain $(\widetilde\rho_n,\widetilde{\bm m}_n,\widetilde S_n)^\top:=\frac{1}{n}\sum_{\ell=1}^n(\rho_\ell,\bm m_\ell,S_\ell)^\top$, where $n=1,\dots,6$.

As explained in \cite{FLY,feireisl2021computing}, a sequence of approximate solutions that converges strongly to a weak solution generates a Young measure concentrated at a single state, that is, a Dirac mass. By contrast, a non-Dirac Young measure can record oscillations that persist under weak convergence. To examine this distinction numerically, we consider a local empirical distribution of the computed density. At a fixed time, this distribution is constructed from the density values $\rho_{j,k}$ at all grid points in the small rectangular subdomain $\widetilde{\Omega}$ introduced below. It therefore measures the spatial variation of the numerical density within this region and provides a common regional diagnostic for the three scheme families. This spatial histogram is different from the pointwise density marginal, which describes the distribution of density states at a fixed space-time point. In particular, the spread of the histogram may reflect both small-scale oscillations and macroscopic spatial variation within $\widetilde{\Omega}$. Its direct identification with the pointwise density marginal would require additional locality assumptions, such as sufficient scale separation or an appropriate shrinking-neighborhood limit. We therefore use the spatial histogram only to compare the regional density distributions produced by the different schemes, rather than as evidence that the limiting Young measure is non-Dirac.

To this end, we take a certain small rectangular subdomain $\widetilde\Omega$ and compute 
$$
\underline\rho:=\min\limits_{(x_j,y_k)\in\widetilde\Omega}\rho_{j,k}\quad\mbox{and}\quad
\xbar\rho:=\max\limits_{(x_j,y_k)\in\widetilde\Omega}\rho_{j,k}.
$$
For each comparison, the same values of $\underline\rho$, $\xbar\rho$, and the same bin edges are used for all distributions being compared. We then introduce $\Delta\rho:=(\xbar\rho-\underline\rho)/30$ and the following partition of the interval $[\underline\rho,\xbar\rho]$:
$$
[\underline\rho,\xbar\rho]=\bigcup\limits_{i=1}^{30}\big[\rho^{(\imh)},\rho^{(\iph)}\big),~~\mbox{where}~~
\rho^{(\hf)}=\underline\rho,~~\rho^{(\iph)}=\rho^{(\imh)}+\Delta\rho,~~i=1,\dots,30.
$$
The distribution of the dataset formed by the values $\rho_{j,k}$ at all grid points $(x_j,y_k)\in\widetilde\Omega$ is represented by a function
$\sigma_\ell=\sigma_\ell(\rho)$, where the subscript $\ell=\ell(r)$ indicates that the density values have been computed using the $r$-th
order scheme. More specifically, $\sigma_\ell(\rho^{(i)})$ denotes the number of grid points $(x_j,y_k)$ for which the computed values $\rho_{j,k}$
fall within the interval $\big[\rho^{(\imh)},\rho^{(\iph)}\big)$. Additionally, we apply a normalization such that 
\begin{equation*}
\Delta\rho\sum\limits_{i=1}^{30}\sigma_\ell\big(\rho^{(i)}\big)=1. 
\end{equation*}
Motivated again by ${\cal K}$-convergence, we average the sequence $\{\sigma_\ell\}$ and introduce
$\widetilde\sigma_n:=\frac{1}{n}\sum\limits_{\ell=1}^n\sigma_\ell$, $n=1,\dots,6$.

We consider the KH instability problem taken from \cite{feireisl2021computing} with the following initial conditions:
$$
(\rho,u,v,p)(x,y,0)=\begin{cases}(2,-0.5,0,2.5),&\mbox{if }0.25+0.05Y_1(x,\omega)<y<0.75+0.05Y_2(x,\omega),\\
(1,0.5,0,2.5),&\mbox{otherwise},\end{cases}
$$
where
$$
Y_j(x,\omega)=\sum_{n=1}^{10}a^n_j(\omega)\cos\big(b^n_j(\omega)+2n\pi x\big),\quad j=1,2,
$$
with $a^n_j=a^n_j(\omega)\in[0,1]$ and $b^n_j=b^n_j(\omega)\in(-\pi,\pi)$ being fixed random numbers, and with $a^n_j$ normalized so that $\sum_{n=1}^{10}a^n_j=1$. The particular values of $a^n_j$ and $b^n_j$ were generated by \texttt{rand} in MATLAB and are reported in \cite[Table 1]{CHKLY26}. The boundary conditions in this example are periodic. In all the Young-measure schemes, we set $[0.2, 3.2]\times [-2,2]\times[-1.5, 1.5]\times[2,10]$ for the discretization in phase space and use 20 grid points in each phase-space variable before retaining only physically admissible nodes. Accordingly, $N_{\bmu}$ denotes the total number of retained phase-space grid points.

\iffalse 
\begin{table}[ht!]
\centering
\begin{adjustbox}{width=\textwidth}
\begin{tabular}{ccccc}
\hline
& $a^n_1$ & $a^n_2$ & $b^n_1$ & $b^n_2$   \\
\hline
1 &6.848086824246653e-08&9.373025805955863e-03&-0.973625473853271& 3.10750325239443\\
2 &4.450348128947341e-03&1.976861219341060e-02& 2.33221742979395 & 1.74829850637860\\
3 &6.156955958786613e-02&1.29928159795144e-01 &-2.57661895600041 &-3.01803367486339\\
4 &1.16481555805349e-01 &2.15403817045303e-01 & 2.43965931651801 &-2.07430785001108\\
5 &1.68204784555961e-01 &1.758678905771403e-02& 1.26278768686501 & 3.10856394704146\\
6 &1.46413246162863e-01 &2.11994809652299e-01 & 1.47373734867445 & 1.59399689987015\\
7 &5.857224323849688e-02&2.134903127162342e-02&-1.25553236484458 &-1.77202310331374\\
8 &1.59868678328304e-01 &1.48283301108284e-01 &-2.82920698380582 & 1.00086476379714\\
9 &1.39003488203723e-01 &1.98161392506709e-01 & 2.56472949845762 & 0.245159802027399\\
10&1.45436027507622e-01 &2.815106156355688e-02&-2.52798558841484 &-0.125265541490568\\
\hline
\end{tabular}
\end{adjustbox}
\caption{\sf The values of $a^n_1$, $a^n_2$, $b^n_1$, and $b^n_2$ for $n=1,\dots,10$ generated by \texttt{rand} in MATLAB.\label{table1a}}
\end{table}
\fi 

\subsection{Numerical Solutions and Approximate Young Measures}
\begin{figure}[ht!]
\centerline{\includegraphics[trim=0cm 0cm 0cm 0cm, clip, width=\linewidth]{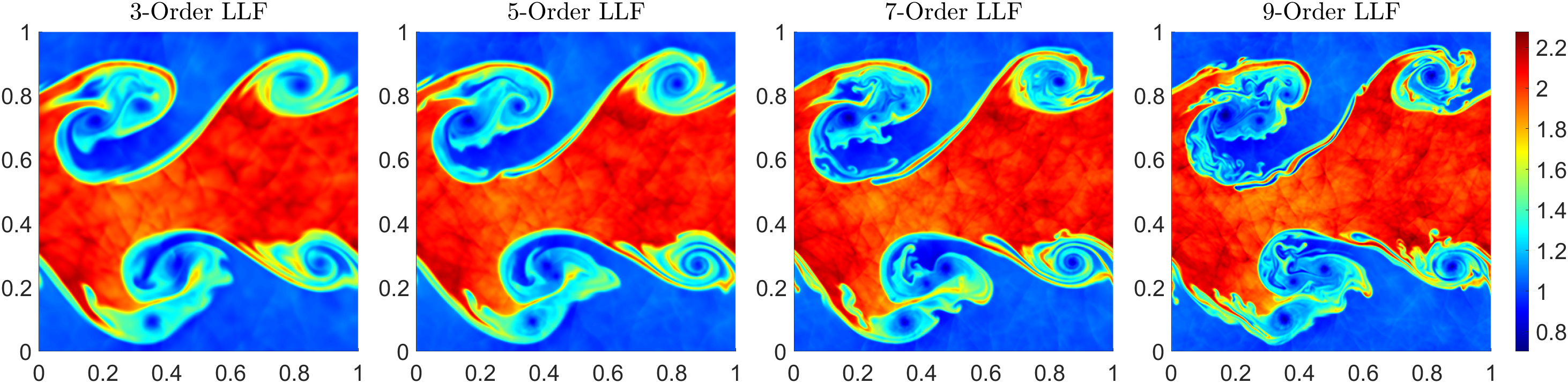}}
\vskip 10pt
\centerline{\includegraphics[trim=0cm 0cm 0cm 0cm, clip, width=\linewidth]{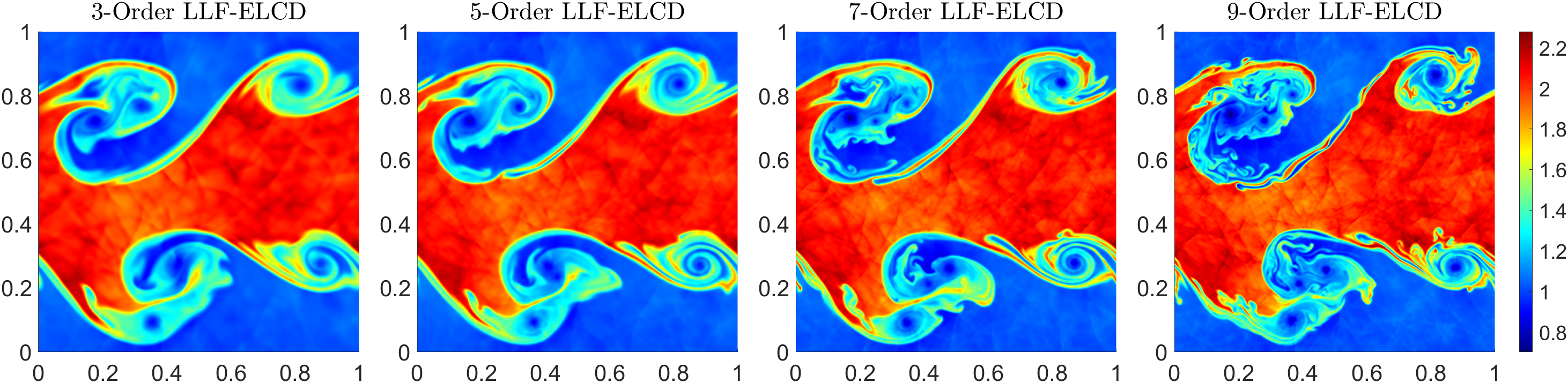}}
\vskip 10pt
\centerline{\includegraphics[trim=0cm 0cm 0cm 0cm, clip, width=\linewidth]{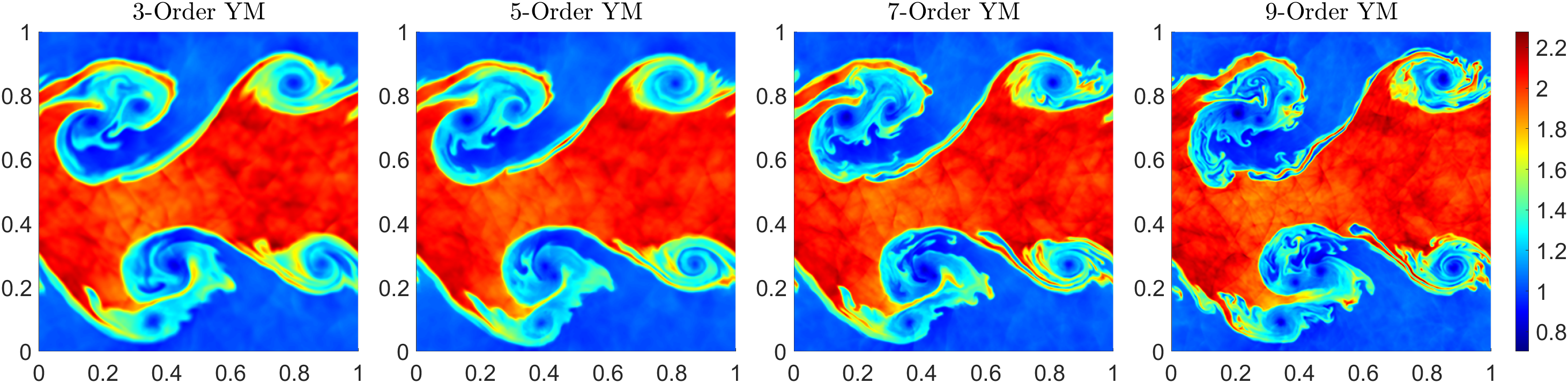}}
\caption{\sf KH Instability: Density computed by the 3-Order (left column), 5-Order (second column), 7-Order (third column), and 9-Order (right column) LLF (top row), LLF-ELCD (middle row), and YM (bottom row) schemes.\label{fig1}}
\end{figure}

We compute the numerical solutions using all the studied schemes until the final time $T=2$ and present in Figure~\ref{fig1} the results obtained by the 3-Order, 5-Order, 7-Order and 9-Order schemes. The results show that the shear layers roll up into vortices, with the higher-order schemes resolving an increasing number of small-scale structures. Moreover, the differences between the LLF and LLF-ELCD schemes are quite limited, while the differences between these two families and the YM schemes are noticeable in the resulting vortical patterns, especially for higher-order schemes. Consequently, at this fixed mesh resolution, the numerical approximations do not appear to stabilize toward a common flow pattern as the formal order increases. Although this finite-order comparison does not establish convergence to distinct DW solutions or rule out strong convergence, it provides evidence of scheme-dependent numerical selection. The cumulative averages and empirical density distributions introduced above will therefore be used to investigate this behavior in greater detail.

\begin{figure}[ht!]
\centerline{\includegraphics[trim=0cm 0cm 0cm 0cm, clip, width=\linewidth]{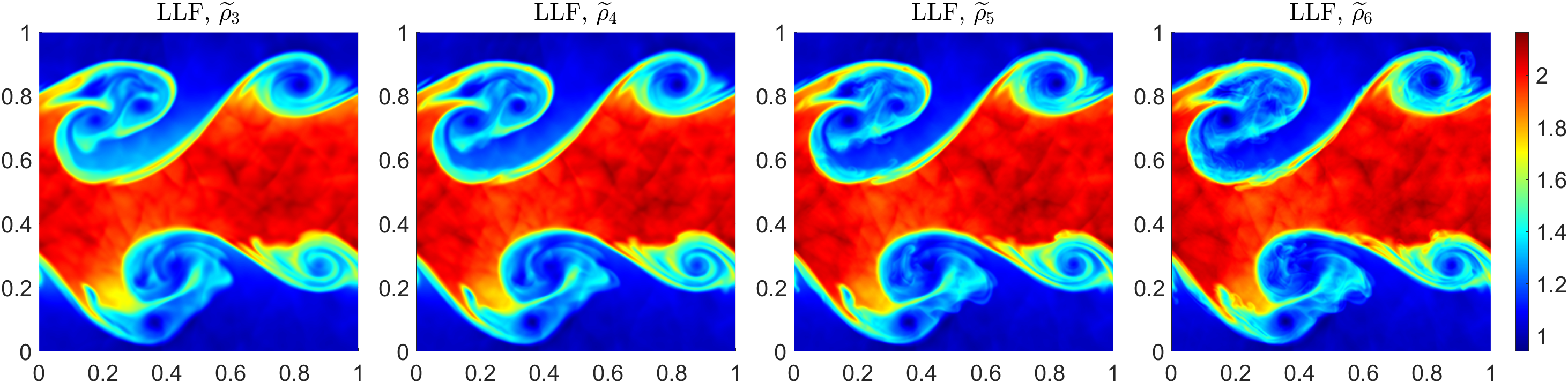}}
\vskip 10pt
\centerline{\includegraphics[trim=0cm 0cm 0cm 0cm, clip, width=\linewidth]{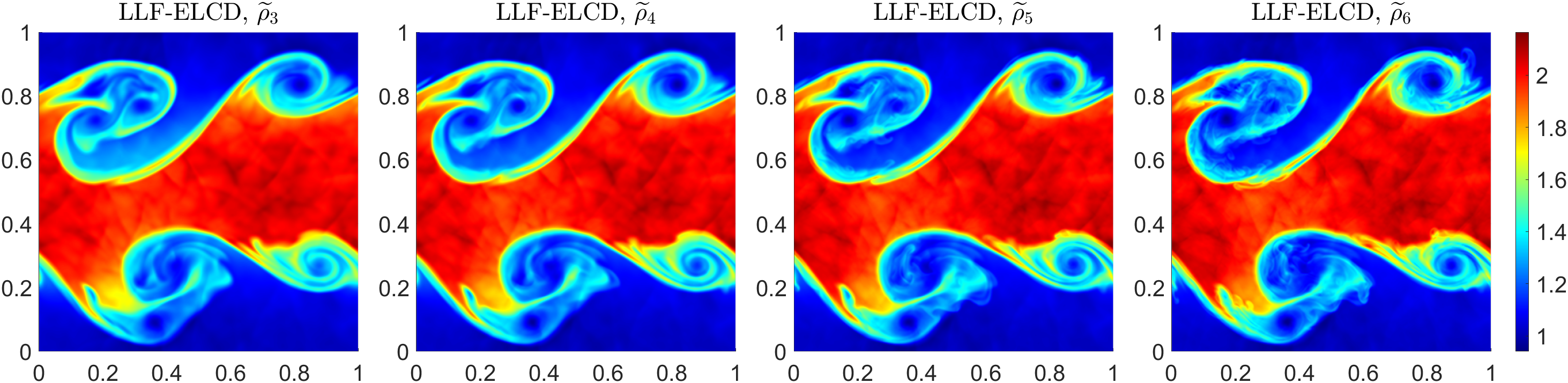}}
\vskip 10pt
\centerline{\includegraphics[trim=0cm 0cm 0cm 0cm, clip, width=\linewidth]{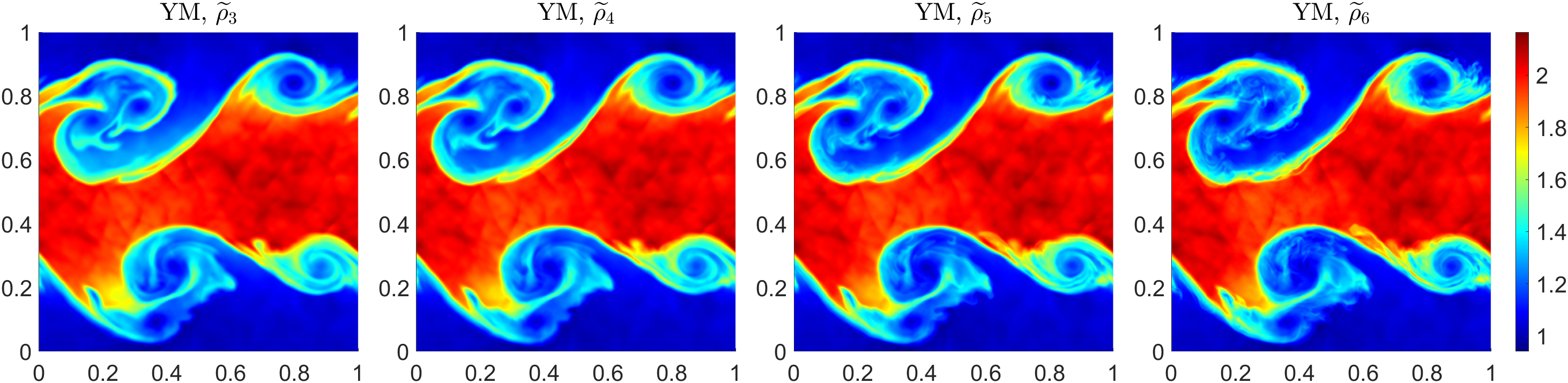}}
\caption{\sf KH Instability: $\widetilde\rho_3$ (left column), $\widetilde\rho_4$ (second column), $\widetilde\rho_5$ (third column), and $\widetilde\rho_6$ (right column) computed by the LLF (top row), LLF-ELCD (middle row), and YM (bottom row)   schemes.\label{fig3}}
\end{figure}

We consider the finite cumulative averages of the numerical solutions over different orders. In Figure \ref{fig3}, we present  $\widetilde\rho_3$, $\widetilde\rho_4$, $\widetilde\rho_5$, and $\widetilde\rho_6$, where one can see that, in contrast to Figure \ref{fig1}, the average density profiles show similarities for different $n$, indicating increased stability as $n$ increases. Additionally, we observe that compared to the LLF and LLF-ELCD schemes, the average density profiles produced by the YM schemes differ significantly. This suggests scheme-dependent numerical selection, although it does not prove convergence to distinct DW solutions.

\begin{figure}[ht!]
\centerline{\includegraphics[trim=0cm 0cm 0cm 0cm, clip, width=\linewidth]{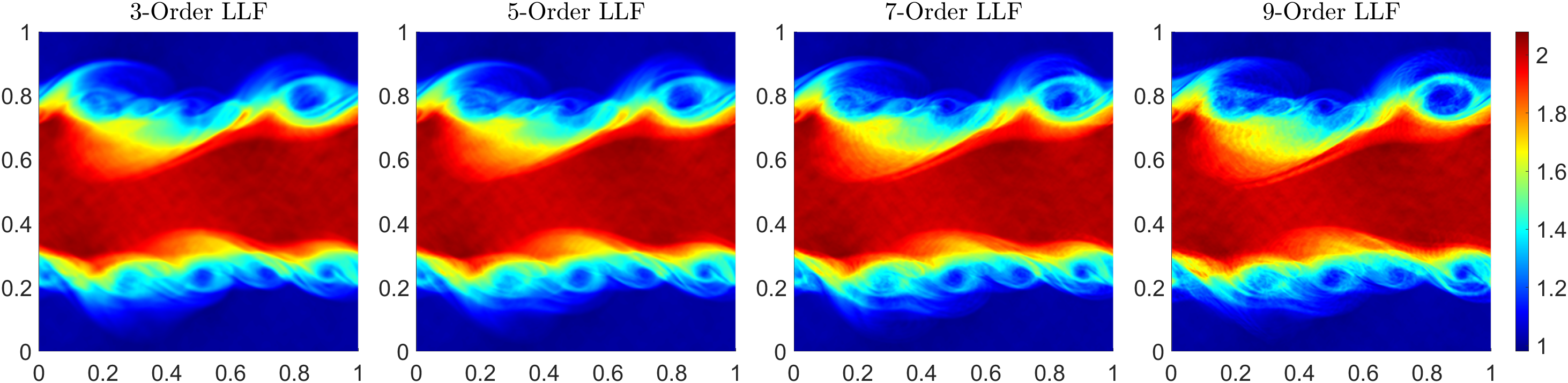}}
\vskip 10pt
\centerline{\includegraphics[trim=0cm 0cm 0cm 0cm, clip, width=\linewidth]{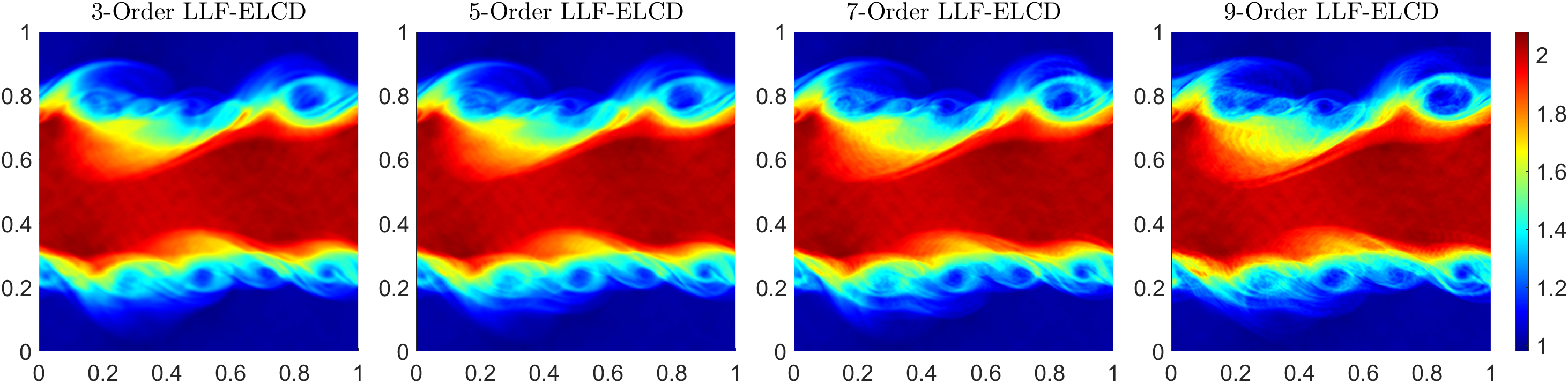}}
\vskip 10pt
\centerline{\includegraphics[trim=0cm 0cm 0cm 0cm, clip, width=\linewidth]{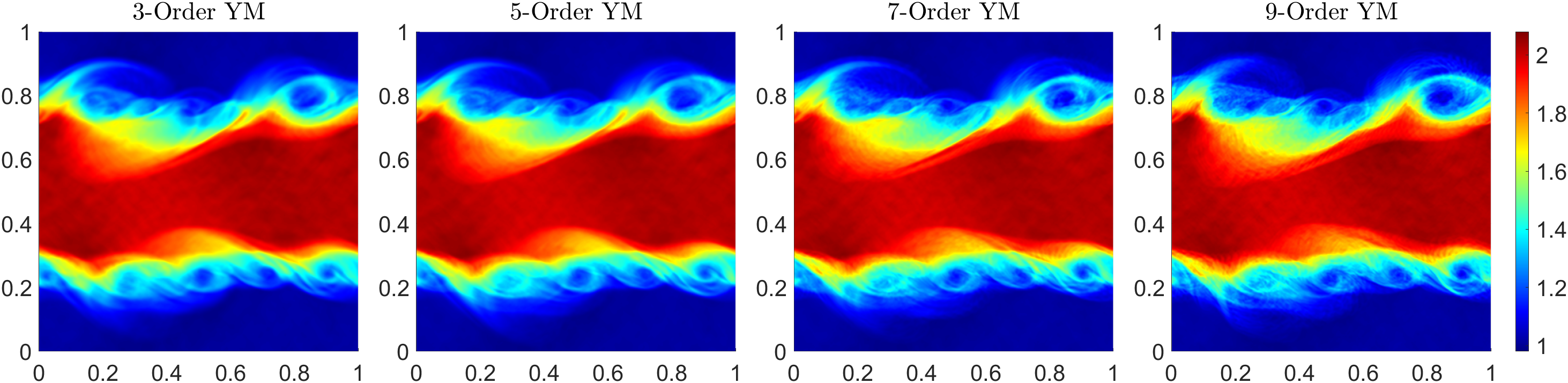}}
\caption{\sf KH Instability: $\rho^T_3$ (left column), $\rho^T_4$ (second column), $\rho^T_5$ (third column), and $\rho^T_6$ (right column) computed by the LLF (top row),
LLF-ELCD (middle row), and YM (bottom row) schemes.\label{fig4}}
\end{figure}

\begin{figure}[ht!]
\centerline{\includegraphics[trim=0cm 0cm 0cm 0cm, clip, width=\linewidth]{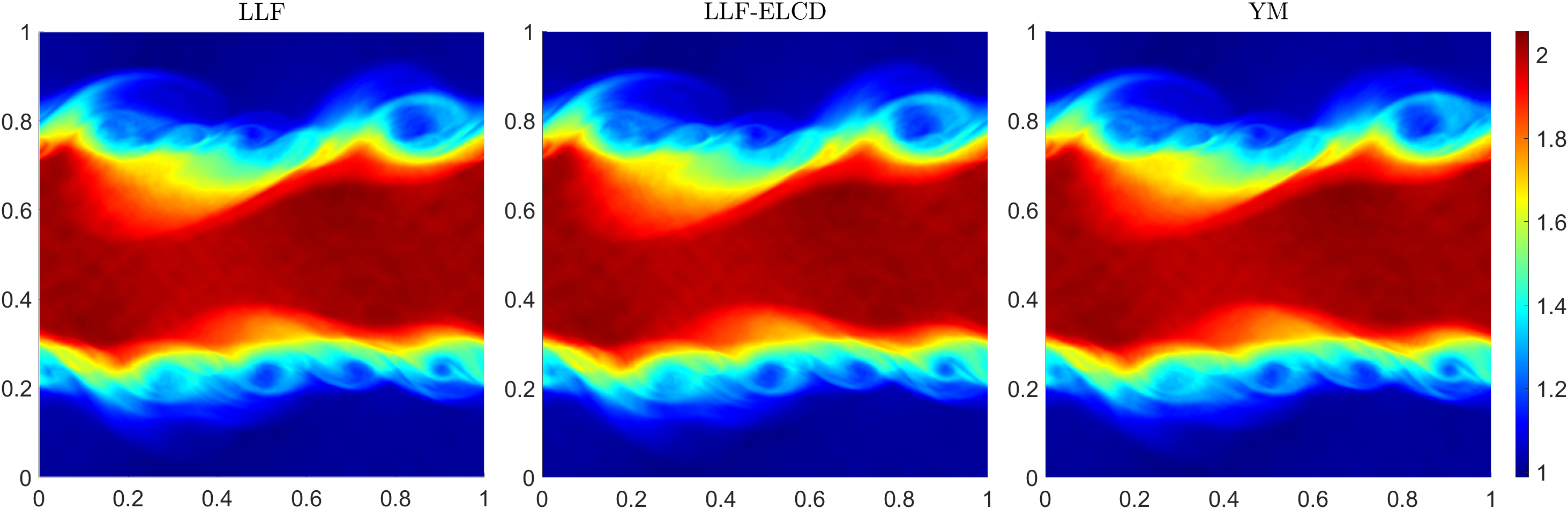}}
\caption{\sf KH Instability: $\widetilde\rho^{\,T}_6$ computed by the LLF (left), LLF-ELCD (middle), and YM (right) schemes.\label{fig5}}
\end{figure}

Next, we approximate the time averages of a DW solution by evaluating the time averages of the computed densities, which are denoted by
$\rho^T_\ell:=\frac{1}{T}\int\limits_0^T\rho_\ell(\cdot,t)\,{\rm d}t$,
where $T$ is the final time. We report $\rho^T_\ell$ for $\ell=3,4,5,6$ computed by the three studied schemes in Figure \ref{fig4}, where one can observe that they become visually more stable as the order increases. In addition, the differences between the time-averaged flows are substantially smaller than those between the instantaneous solutions. 
In Figure \ref{fig5}, we show the time averages of $\widetilde{\rho}_6$, which are denoted by $\widetilde\rho^{\,T}_6:=\frac{1}{T}\int\limits_0^T\widetilde\rho_6(\cdot,t)\,{\rm d}t$. The three time-averaged profiles are remarkably similar. This phenomenon suggests consistency in the measure-valued solutions produced by different numerical methods, but does not establish that the three methods generate the same measure-valued limit.

\begin{figure}[ht!]
\centerline{\includegraphics[trim=1.1cm 0.3cm 1.4cm 0.1cm, clip, width=4.8cm]{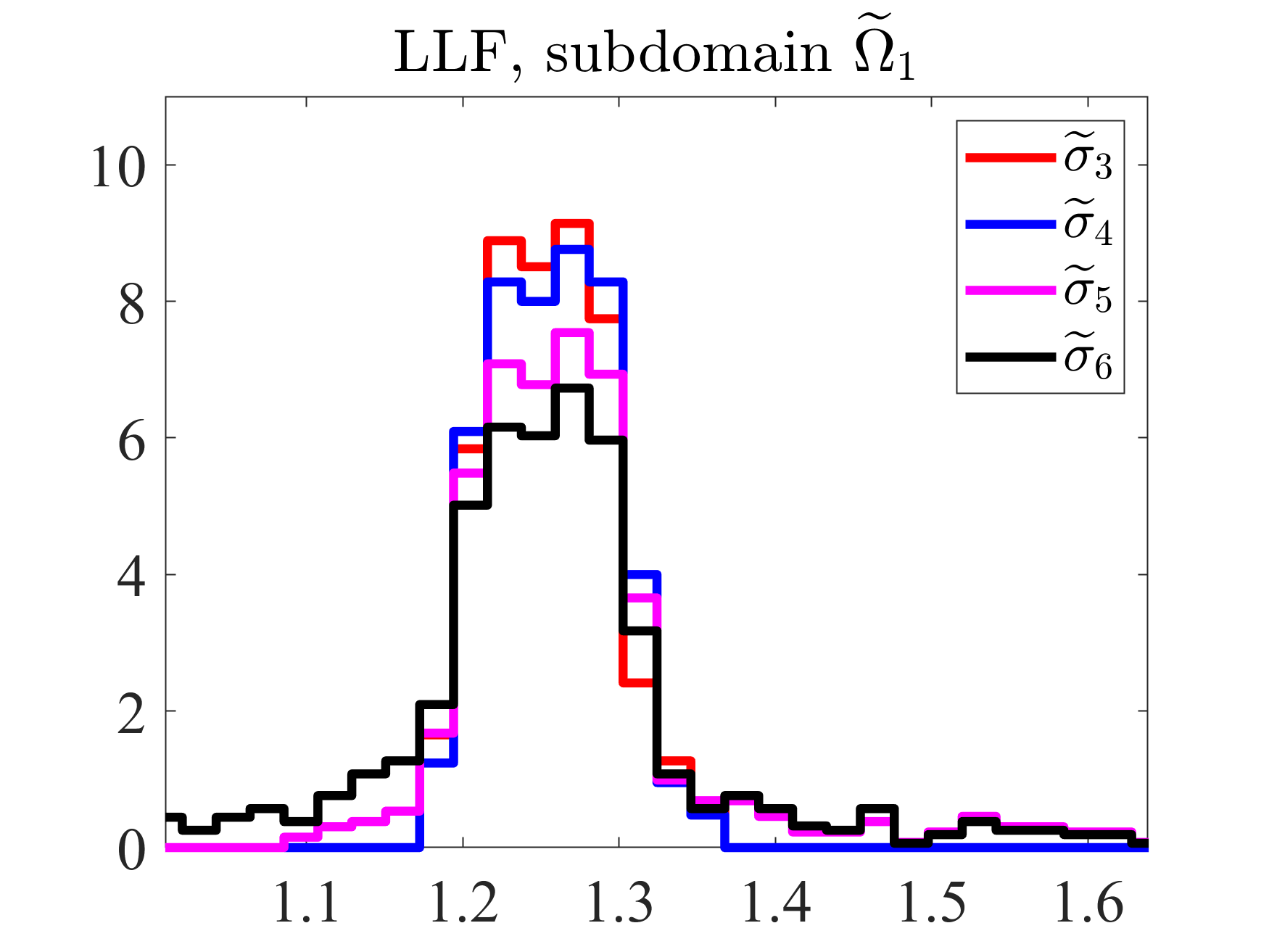}\hspace*{0.5cm}
            \includegraphics[trim=1.1cm 0.3cm 1.4cm 0.1cm, clip, width=4.8cm]{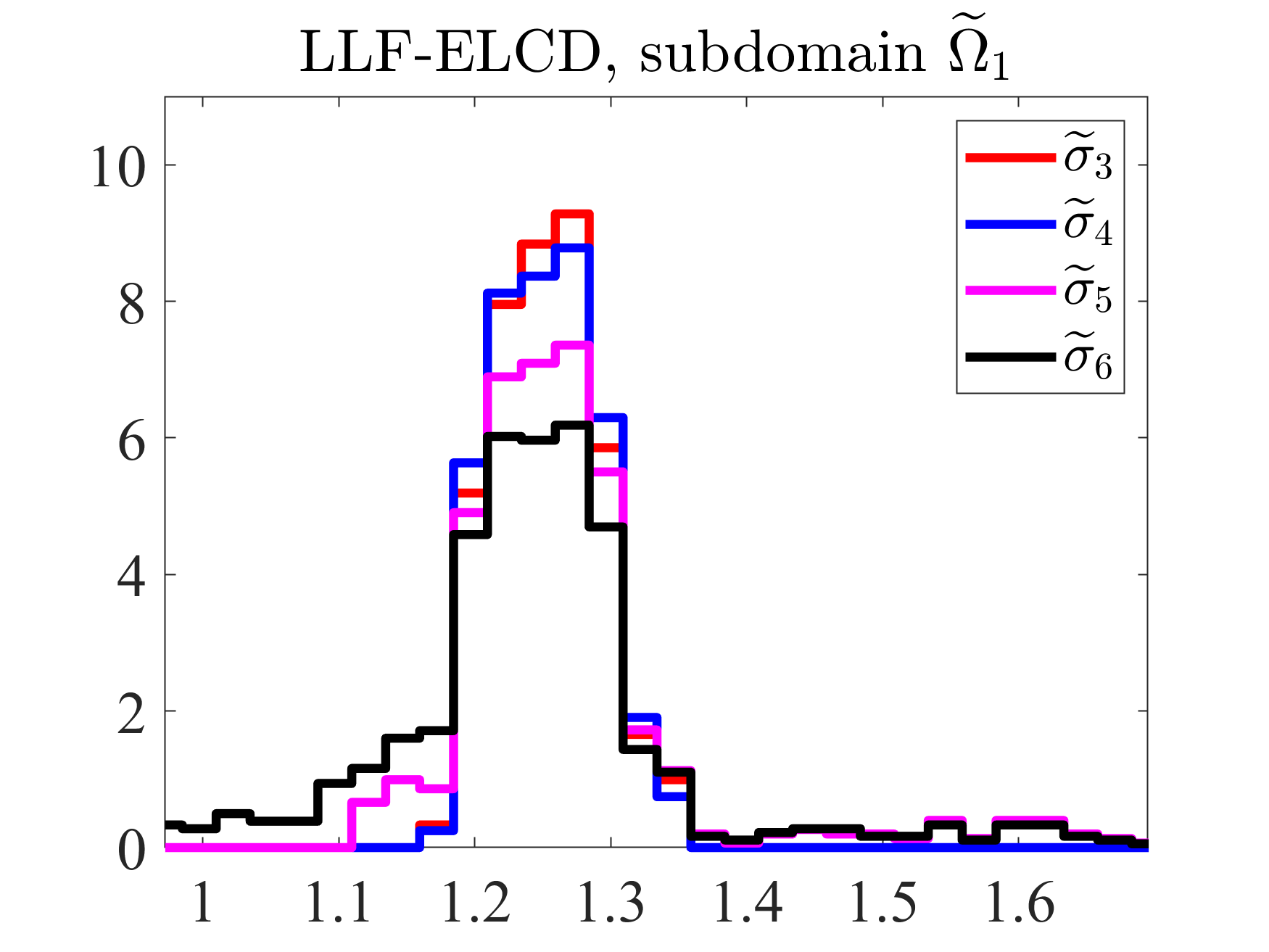}\hspace*{0.5cm}
            \includegraphics[trim=1.1cm 0.3cm 1.4cm 0.1cm, clip, width=4.8cm]{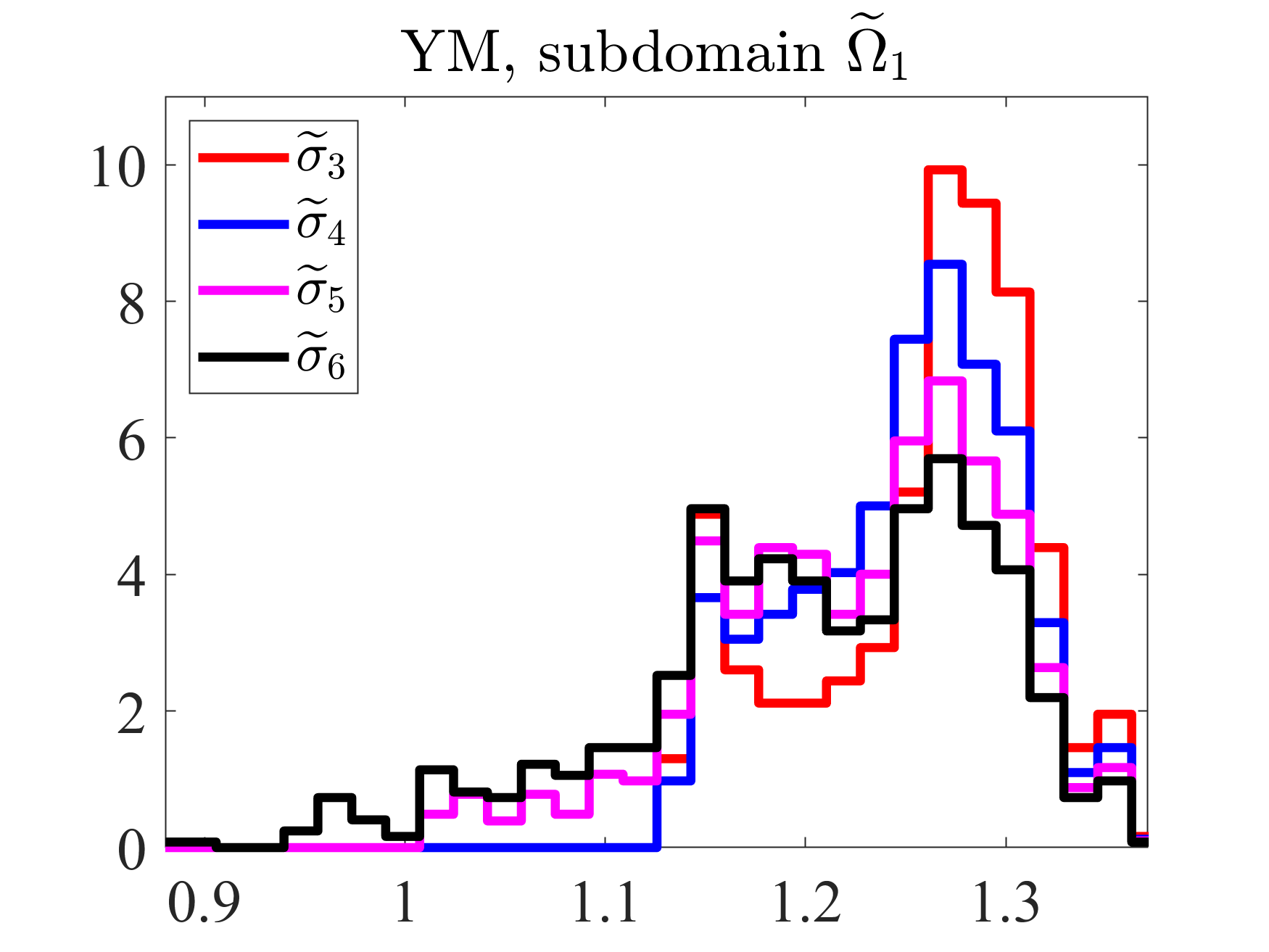}}
\vskip 10pt
\centerline{\includegraphics[trim=1.1cm 0.3cm 1.4cm 0.1cm, clip, width=4.8cm]{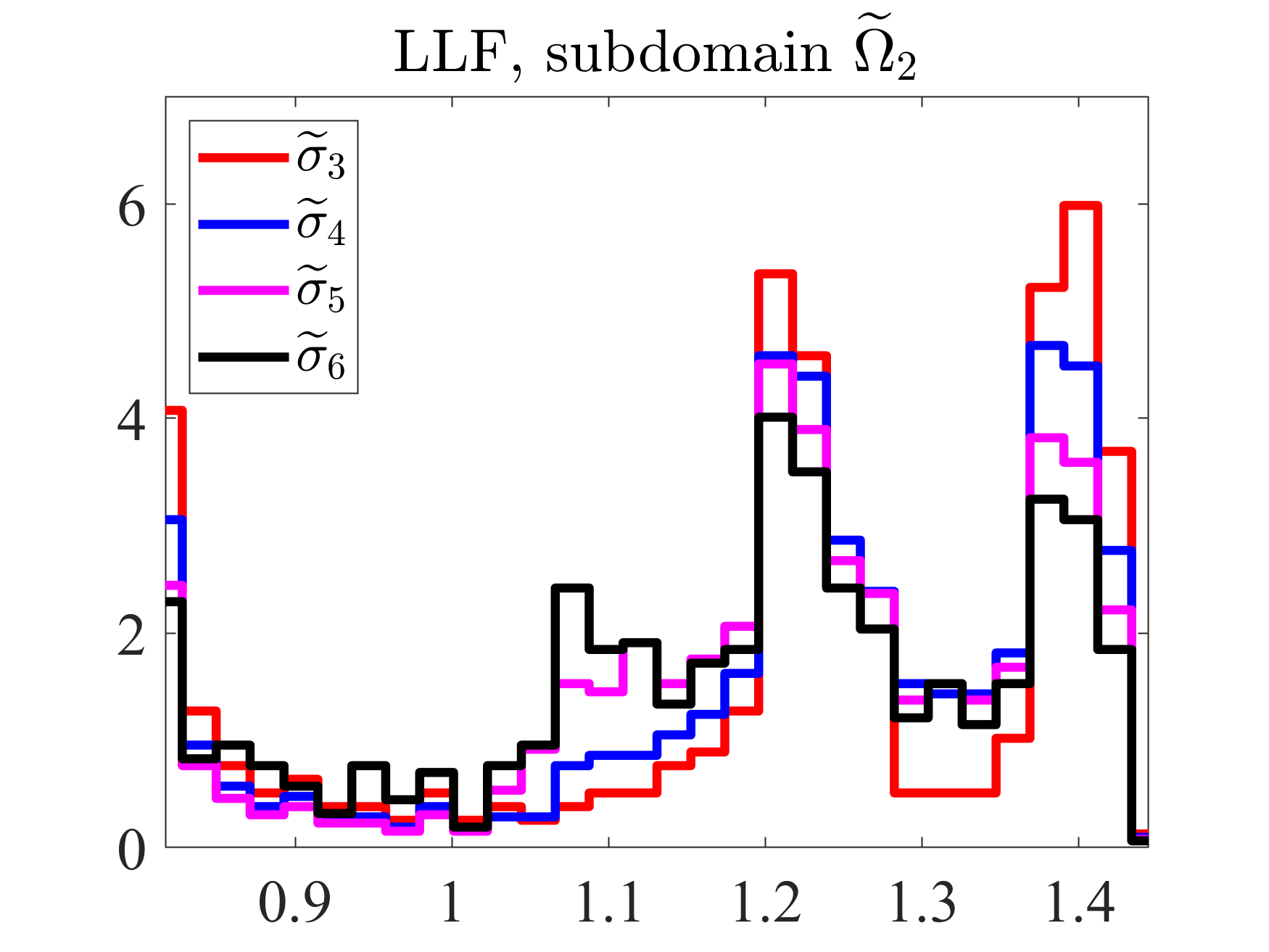}\hspace*{0.5cm}
            \includegraphics[trim=1.1cm 0.3cm 1.4cm 0.1cm, clip, width=4.8cm]{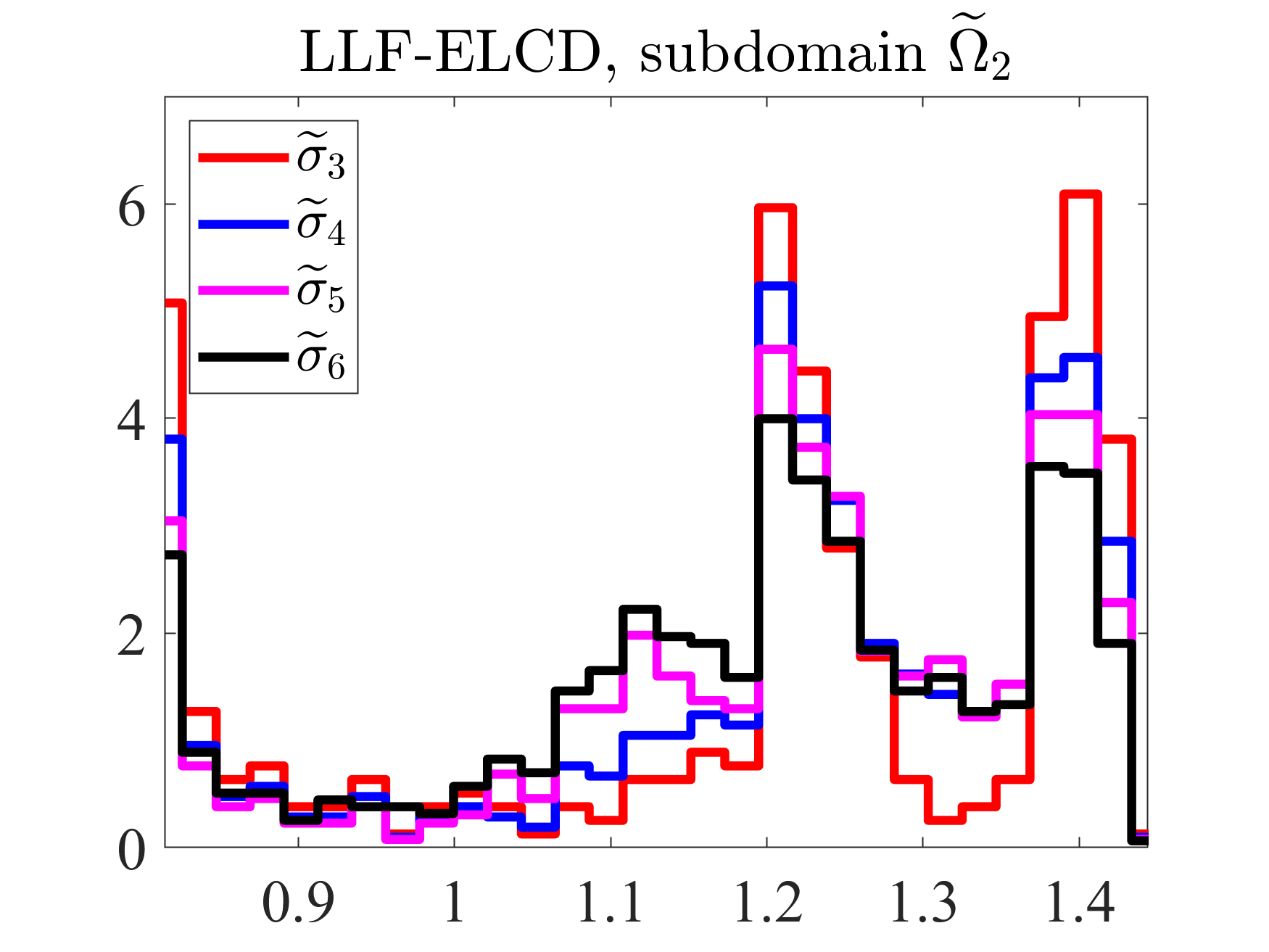}\hspace*{0.5cm}
            \includegraphics[trim=1.1cm 0.3cm 1.4cm 0.1cm, clip, width=4.8cm]{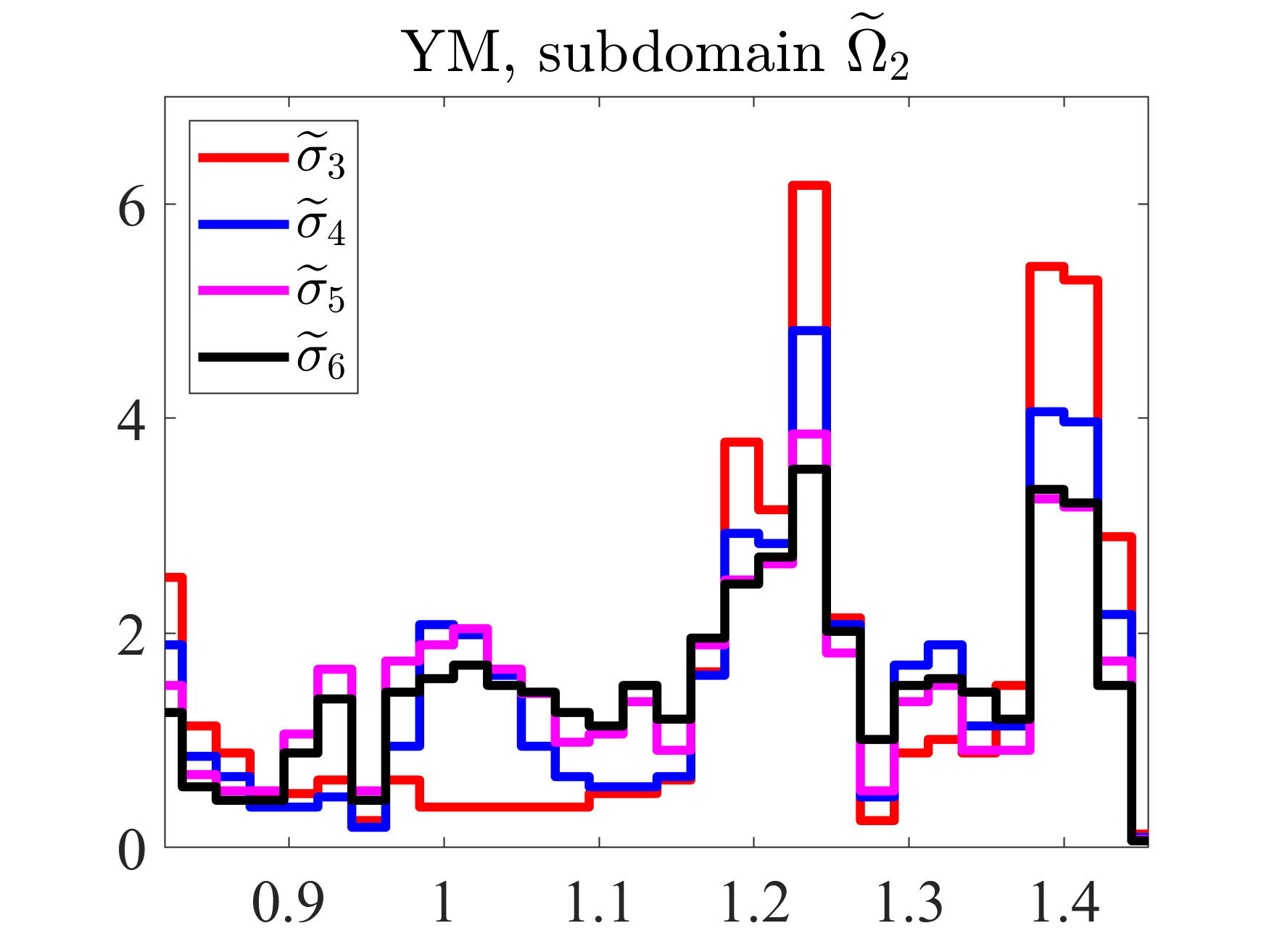}}
\caption{\sf KH Instability: $\widetilde\sigma_3$, $\widetilde\sigma_4$, $\widetilde\sigma_5$, and $\widetilde\sigma_6$ computed by the
LLF (left column), LLF-ELCD (middle column), and YM (right column) schemes in the subdomains $\widetilde\Omega_1$ (top row) and
$\widetilde\Omega_2$ (bottom row).\label{fig6}}
\end{figure}

To illustrate the differences in numerical results, we also compare the regional probability density functions (PDFs) computed by the three studied schemes. To this end, we select two subdomains:
\begin{equation*}
\widetilde\Omega_1=[0.24,0.26]\times[0.80,0.82]\quad\mbox{and}\quad\widetilde\Omega_2=[0.36,0.38]\times[0.24,0.26],
\end{equation*}
and present the obtained results in Figure \ref{fig6}. It is evident that these PDFs computed by the LLF and LLF-ELCD schemes have similarities, especially in the subdomain $\widetilde\Omega_2$, while those computed by the YM method differ significantly from the other two. This provides additional evidence that different methods may select different DW solutions. 

\subsection{Comparison Based on Selection Criteria}
According to \cite{FL26,FLY,CHKLY26}, possible selection criteria for DW solutions are (i) the maximization of entropy production, (ii) the minimization of the energy defect, (iii) the maximization of the energy defect, and (iv) the minimization of a Bregman functional associated with a prescribed constant reference state. The first criterion favors larger physical entropy, while the fourth combines the mean-state energy and physical entropy through the functional defined below. The second and third criteria express different preferences for the energy defect, which controls the Reynolds stress defect in the DW formulation.  We refer to \cite{Lukacova_book,FL26} for the precise relation between these defects and the corresponding selection results.  In this paper, we use these criteria to compare the numerical solutions obtained by the three scheme families. For this purpose, we compute the time-averaged total physical entropy, an energy defect obtained from the finite cross-order averages, and a Bregman functional. These quantities provide numerical comparisons at the present resolution; they do not by themselves determine a selected DW limit.

\begin{figure}[ht!]
\centerline{\includegraphics[trim=0cm 0cm 0cm 0cm, clip, width=\linewidth]{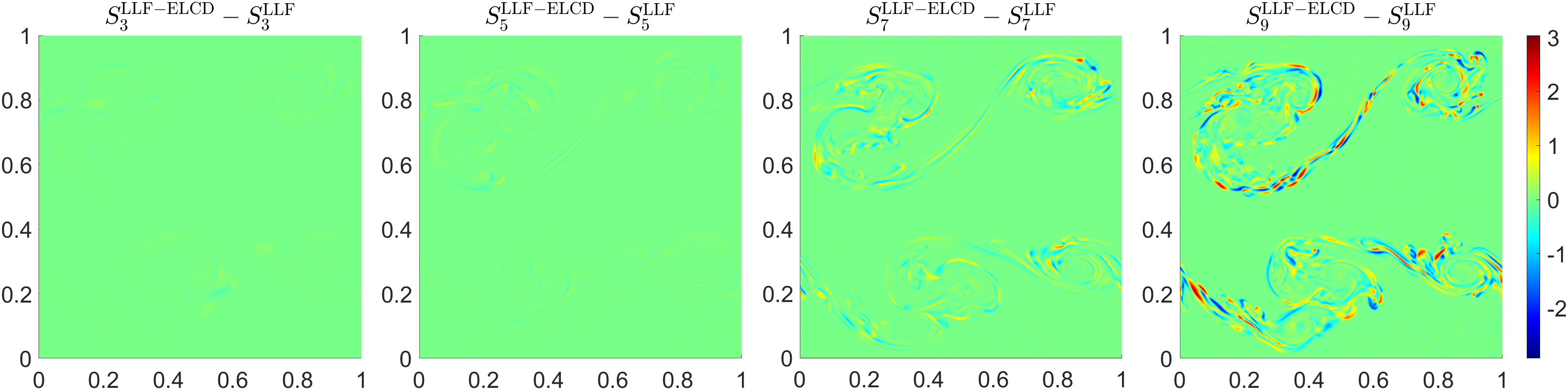}}
\vskip 10pt
\centerline{\includegraphics[trim=0cm 0cm 0cm 0cm, clip, width=\linewidth]{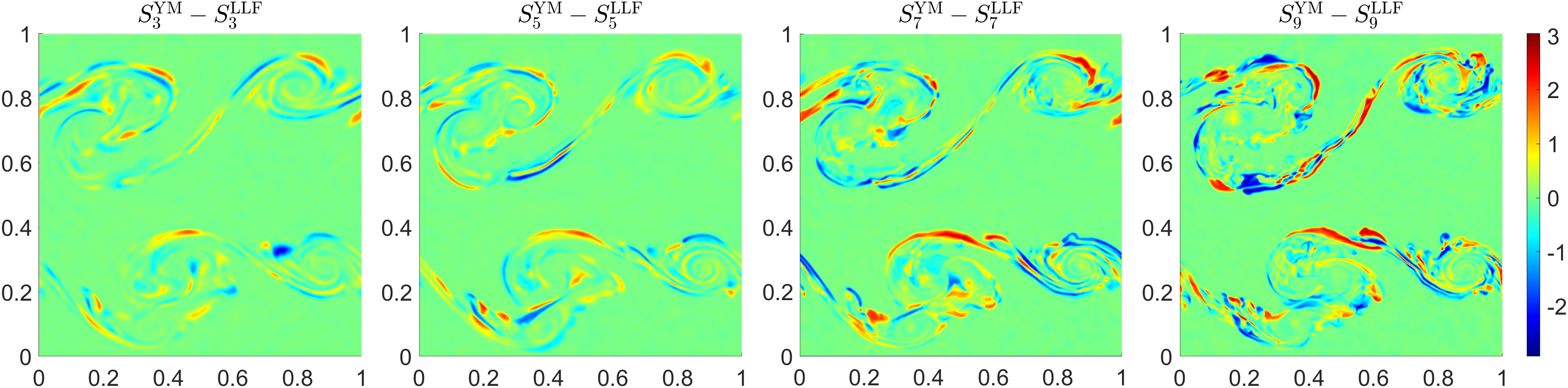}}
\caption{\sf KH Instability: The entropy differences between the 3-Order (left column), 5-Order (second column), 7-Order (third column), and 9-Order (right column) LLF with the corresponding LLF-ELCD (top row) and YM (bottom row) schemes at the final time $T=2$.\label{fig44}}
\end{figure}

We first examine the entropy-related selection criterion. In the LLF-ELCD schemes, the entropy-based LCD uses physical entropy only to select the interface state that defines the local characteristic basis; see Appendix \ref{appb}. In the YM schemes, the expected mathematical entropy is minimized in the LP problem \eref{linprog-a}--\eref{linprog-e} to compute the discrete probability measures. We therefore examine whether these two uses of entropy lead to larger physical entropy in the computed solutions. In Figure ~\ref{fig44}, we show the differences in the entropy at the final time among the LLF, LLF–ELCD, and YM schemes. We observe that these differences are generally more pronounced as the formal order increases. For the lower-order schemes, the differences between the LLF and LLF-ELCD results are relatively small. By contrast, the differences between the LLF and YM results are substantially more pronounced than those between the LLF and LLF-ELCD results.

In Table \ref{table1}, we report the total physical entropy for each formal order $r$,
\begin{equation*}
{\cal S}_r(t):=\int\limits_\Omega S_r(x,y,t)\,{\rm d}x{\rm d}y,
\end{equation*}
computed by the three schemes studied. Here, the LLF-ELCD schemes do not achieve larger average entropies of different formal orders, while the YM schemes produce larger entropies in {\bf all} of the formal orders. The numerical results demonstrate that using entropy only to select an interface state for the LCD does not optimize the numerical solution and also does not enforce a maximum-entropy selection principle. By contrast, the YM results are consistent with the local maximum-entropy bias imposed through the objective function in the LP problem, but they do not by themselves prove selection of a globally maximal-entropy DW solution.
\begin{table}[ht!]
\centering
\caption{\sf KH Instability: Time averages of the total physical entropy for the 1-Order, 2-Order, 3-Order, 5-Order, 7-Order, and 9-Order LLF, LLF-ELCD, and YM schemes.\label{table1}}
\begin{tabular}{ccccc}
\toprule 
& \multirow{2}{*}{LLF} &\multirow{2}{*}{LLF--ELCD}& \multirow{2}{*}{YM} & Relative improvement of\\& & & & YM over LLF (\%)\\ \hline 
$\frac{1}{T}\int_0^T{\cal S}_1(t)\,{\rm d}t$ &1.1160477    & 1.1160477   &1.1177841 & 0.15558\\\hline
$\frac{1}{T}\int_0^T{\cal S}_2(t)\,{\rm d}t$ &1.0627867    & 1.0627606   &1.0632212 & 0.04088\\\hline
$\frac{1}{T}\int_0^T{\cal S}_3(t)\,{\rm d}t$ &1.0617010    & 1.0616969   &1.0623086 & 0.05723\\\hline
$\frac{1}{T}\int_0^T{\cal S}_5(t)\,{\rm d}t$ &1.0609902    & 1.0609699   &1.0629257 & 0.18242\\\hline
$\frac{1}{T}\int_0^T{\cal S}_7(t)\,{\rm d}t$ &1.0556312    & 1.0555922   &1.0582666 & 0.24965\\\hline
$\frac{1}{T}\int_0^T{\cal S}_9(t)\,{\rm d}t$ &1.0528719    & 1.0527596   &1.0543940 & 0.14457\\
\bottomrule 
\end{tabular}
\end{table}

We then compute the average total entropy over orders up to $n$ by
\begin{equation*}
{\cal S}^n(t):=\int\limits_\Omega\widetilde S_n(x,y,t)\,{\rm d}x{\rm d}y,
\end{equation*}
the averaged total energy and the energy of the cumulative mean state by
\begingroup
\begin{equation*}
\widetilde E_n(x,y,t):=\frac{1}{n}\sum_{\ell=1}^nE\big(\rho_\ell(x,y,t),\bm m_\ell(x,y,t),S_\ell(x,y,t)\big),
\end{equation*}
\endgroup
\begin{equation*}
\begin{aligned}
{\cal E}^n_1(t)&:=\int\limits_\Omega\widetilde E_n(x,y,t)\,{\rm d}x{\rm d}y=
\int\limits_\Omega\frac{1}{n}\sum_{\ell=1}^n\left[\frac{|\bm m_\ell|^2}{2\rho_\ell}+\rho_\ell e(\rho_\ell,S_\ell)\right]{\rm d}x{\rm d}y,\\
{\cal E}^n_2(t)&:=\int\limits_\Omega E\big(\widetilde\rho_n(x,y,t),\widetilde{\bm m}_n(x,y,t),\widetilde S_n(x,y,t)\big){\rm d}x{\rm d}y=
\int\limits_\Omega\frac{|\widetilde{\bm m}_n|^2}{2\widetilde\rho_n}+\widetilde\rho_ne\big(\widetilde\rho_n,\widetilde S_n\big)
{\rm d}x{\rm d}y,
\end{aligned}
\end{equation*}
respectively, and the energy defect is defined by
\begin{equation*}
{\cal D}^n_{\rm E}(t):= \int\limits_\Omega\left|\widetilde E_n(x,y,t)-E\big(\widetilde\rho_n(x,y,t),\widetilde{\bm m}_n(x,y,t),
\widetilde S_n(x,y,t)\big)\right|{\rm d}x{\rm d}y.
\end{equation*}
In addition, we consider a Bregman functional associated with the constant reference state $(\bar\rho,\bm 0,S(\bar\rho,\bar\theta))$ that maximizes the entropy. Here,
$$
\bar{\rho}= \frac{1}{|\Omega|} \int\limits_\Omega \rho_0 \, \mathrm{d}x\mathrm{d}y,
\quad \bar\theta=\frac{1}{ c_v \bar\rho |\Omega|} \int\limits_\Omega E_0 \, \mathrm{d}x\mathrm{d}y,
$$
and the Bregman functional is given by
\begin{equation*}
{\mathcal J}^n(t):= \mathcal{E}^n_2(t) -\bar \theta \mathcal{S}^n(t);
\end{equation*}
see \cite[\S6]{FL26} for details on its derivation.

Based on the numerical solutions obtained by different methods, we compare their respective time-averaged total entropy, energy defects, and the Bregman functional. The time averages of ${\cal S}^6$, ${\mathcal E}_1^6$, ${\mathcal E}_2^6$, ${\mathcal D}^6_{\rm E}$, and $\mathcal{J}^6$ are presented in Table \ref{table2}. Here, the initial data yield $\bar{\theta}=1.71665$. It can be observed that the YM schemes produce a smaller energy defect than the LLF and LLF-ELCD schemes. At the same time, these two methods result in smaller time-averaged total entropies.  According to these numerical diagnostics, the YM schemes are prioritized by criteria (i), (ii), and (iv), while the LLF-ELCD schemes are prioritized by criterion (iii). Note that the differences between the LLF and LLF-ELCD schemes with respect to all selection criteria are very small. 

\begin{table}[ht!]
\centering
\caption{\sf KH Instability: Time averages of the total physical entropy, averaged total energy, energy of the cumulative mean state, the energy defect, and the Bregman functional computed by the LLF, LLF-ELCD, and YM schemes.\label{table2}}
\begin{tabular}{ccccc}\toprule 
& \multirow{2}{*}{LLF} &\multirow{2}{*}{LLF--ELCD}& \multirow{2}{*}{YM} & Relative change of\\& & & & YM over LLF (\%)\\ \hline 
$\frac{1}{T}\int_0^T{\cal S}^6(t)\,{\rm d}t$       &1.0683381     &1.0683045     & 1.0698167  & 0.13840\\\hline
$\frac{1}{T}\int_0^T{\mathcal E}^6_1(t)\,{\rm d}t$ &6.4375014     &6.4375014     & 6.4375014  & 0\\\hline
$\frac{1}{T}\int_0^T{\mathcal E}^6_2(t)\,{\rm d}t$ &6.4127463     &6.4127242     & 6.4136133  & 0.01352\\\hline
$\frac{1}{T}\int_0^T{\mathcal D}^6_E(t)\,{\rm d}t$ &2.47551e-2    &2.47772e-2    & 2.38881e-2 & -3.50231\\\hline
$\frac{1}{T}\int_0^T{\mathcal J}^6(t)\,{\rm d}t$   &4.5787794     &4.5788151     & 4.5771083  & -0.03650\\
\bottomrule 
\end{tabular}
\end{table}

\subsection{Density Marginals of the Averaged Young Measures}
At each physical grid point $(x_j,y_k)$, the coefficients $\mu^*_{j,k,\ell}$ describe the local distribution of the Young measure over the discrete phase-space states $\bm z_\ell$. They satisfy
$$
\mu^*_{j,k,\ell}\geq 0,\quad \Delta V\sum_{\ell=1}^{N_{\bm u}}\mu^*_{j,k,\ell}=1,
$$
where $\Delta V$ is the constant phase-space cell volume. Therefore, $\mu^*_{j,k,\ell}$ is a weight of the probability density. We define
$w_{j,k,\ell}:=\Delta V\mu^*_{j,k,\ell}$.  The normalized coefficients $w_{j,k,\ell}$ satisfy $w_{j,k,\ell}\geq0$ and $\sum_{\ell=1}^{N_{\bm u}}w_{j,k,\ell}=1$. Accordingly, the discrete Young measure at $(x_j,y_k)$ can be written as $\nu^*_{j,k}=\sum_{\ell=1}^{N_{\bm u}}w_{j,k,\ell}\delta_{\bm z_\ell}$. Moreover, its first moment satisfies $\sum_{\ell=1}^{N_{\bm u}}w_{j,k,\ell}\bm z_\ell =\bm U_{j,k}$.

The values of $w_{j,k,\ell}$ provide direct information about the local concentration of the Young measure. If a weight is close to one, the measure is strongly concentrated at a single phase-space state. It is close to a Dirac measure. If several weights have comparable nonzero values, the measure is a locally dispersed Young measure.  Therefore, the spatial distributions of $w_{j,k,\ell}$, together with quantities such as $\max_\ell w_{j,k,\ell}$ and the effective number of active states, can be used to examine the concentration of the Young measure.

In the following, we study the density marginals averaged first over a selected physical subdomain and then over several YM reconstruction orders. Thus, the two figures do
not display the local full-state weights $w_{j,k,\ell}$. For the YM scheme of order $r\in\{1,2,3,5,7,9\}$, we write the density marginal at $(x_j,y_k)$ as $\omega_{j,k}^{(r)}=\sum_{m=1}^{N_\rho}w_{j,k,m}^{\rho,(r)}\delta_{\rho_m^*}$,
where $\rho_m^*$ is the $m$-th density support point. The mass $w_{j,k,m}^{\rho,(r)}$ is obtained by summing the weights of all phase-space states whose density component equals $\rho_m^*$. In the numerical output, it is recovered from the saved density-marginal coefficient through $w_{j,k,m}^{\rho,(r)}=\Delta\rho\,\mu_{j,k,m}^{\rho,(r)}$,
where $\mu_{j,k,m}^{\rho,(r)}$ is the density-marginal coefficient. These marginal masses satisfy
\begin{equation*}
 w_{j,k,m}^{\rho,(r)}\geq0,\quad \sum_{m=1}^{N_\rho}w_{j,k,m}^{\rho,(r)}=1,\quad\sum_{m=1}^{N_\rho}\rho_m^*w_{j,k,m}^{\rho,(r)}=\rho_{j,k}^{(r)}.
\end{equation*}

We next perform two successive averages: First, we average over physical cells in a selected spatial subdomain. Let $\mathcal I$ denote the corresponding set of grid indices corresponding to the physical cells. For the scheme of order $r$, we define
\begin{equation*}
 \overline w_{m}^{\rho,(r)}:=\frac{1}{|\mathcal I|}\sum_{(j,k)\in\mathcal I}w_{j,k,m}^{\rho,(r)},\quad m=1,\ldots,N_\rho.
\end{equation*}
The same definition is applied separately to each selected subdomain. Since the local density marginals are normalized, their spatial averages satisfy
$$
\overline w_{m}^{\rho,(r)}\geq0,\quad \sum_{m=1}^{N_\rho}\overline w_{m}^{\rho,(r)}=1.
$$
We then average the spatially averaged masses obtained with different YM reconstruction orders. To this end, we introduce
$$
\mathcal R_3=\{1,2,3\},\quad 
\mathcal R_4=\{1,2,3,5\},\quad 
\mathcal R_5=\{1,2,3,5,7\},\quad {\rm and} \quad 
\mathcal R_6=\{1,2,3,5,7,9\}.
$$
For $q=3,\ldots,6$, we define
\begin{equation*}
 \widetilde w_{q,m}^{\rho}:=\frac{1}{q}\sum_{r\in\mathcal R_q}\overline w_{m}^{\rho,(r)},\quad
 \widetilde\omega_{q}:=\sum_{m=1}^{N_\rho}\widetilde w_{q,m}^{\rho}\delta_{\rho_m^*}.
\end{equation*}
Here, $q$ denotes the number of schemes included in the average and not the formal order of a numerical scheme. Note that every $\widetilde\omega_{q}$ is a probability measure on $\mathbb{R}$.

\begin{figure}[ht!]
\centerline{\includegraphics[trim=0.0cm 0.cm  1.4cm 0.1cm, clip, width=7.5cm]{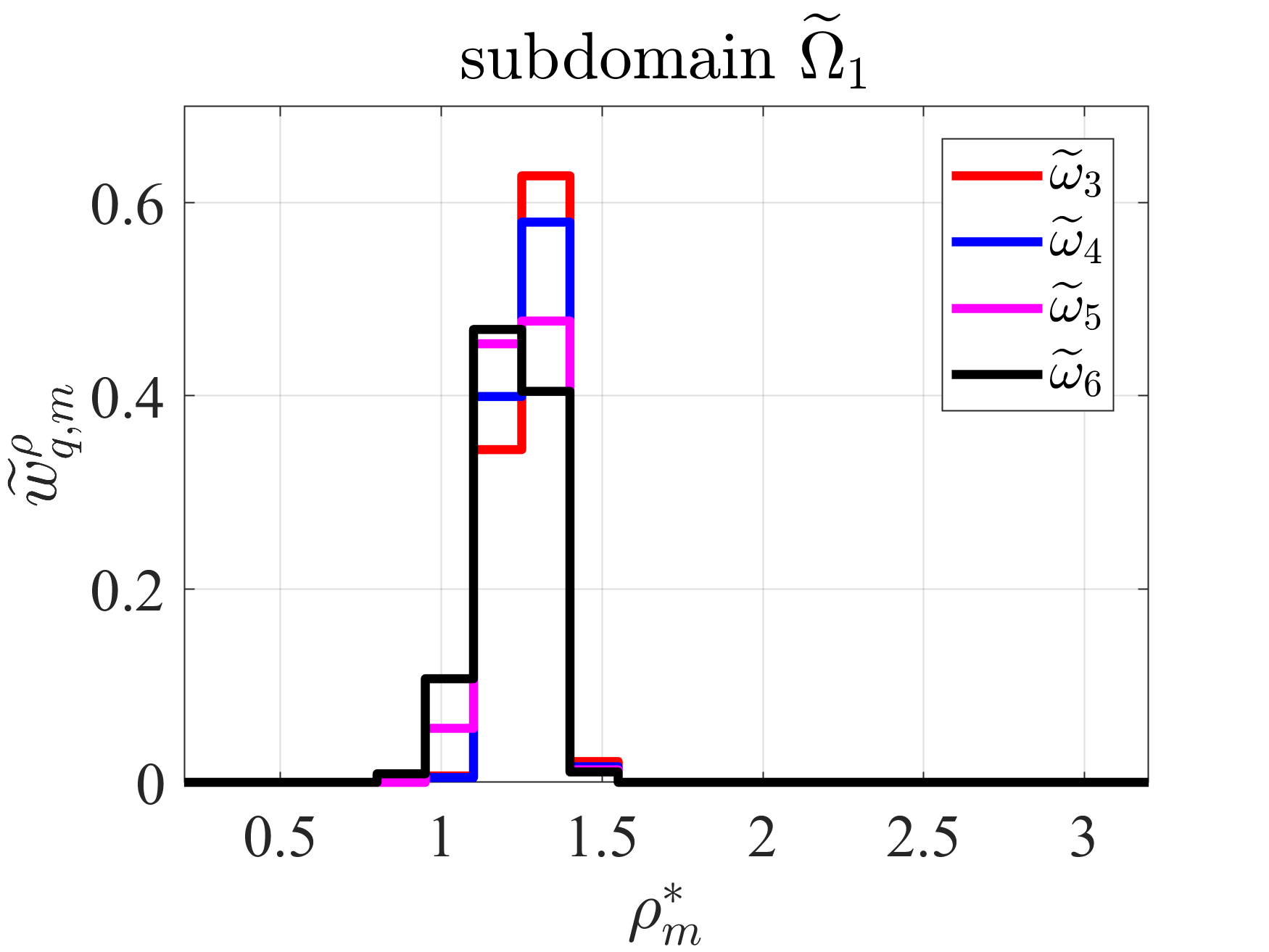}\hspace*{1.cm}
            \includegraphics[trim=0.0cm 0.0cm 1.4cm 0.1cm, clip, width=7.5cm]{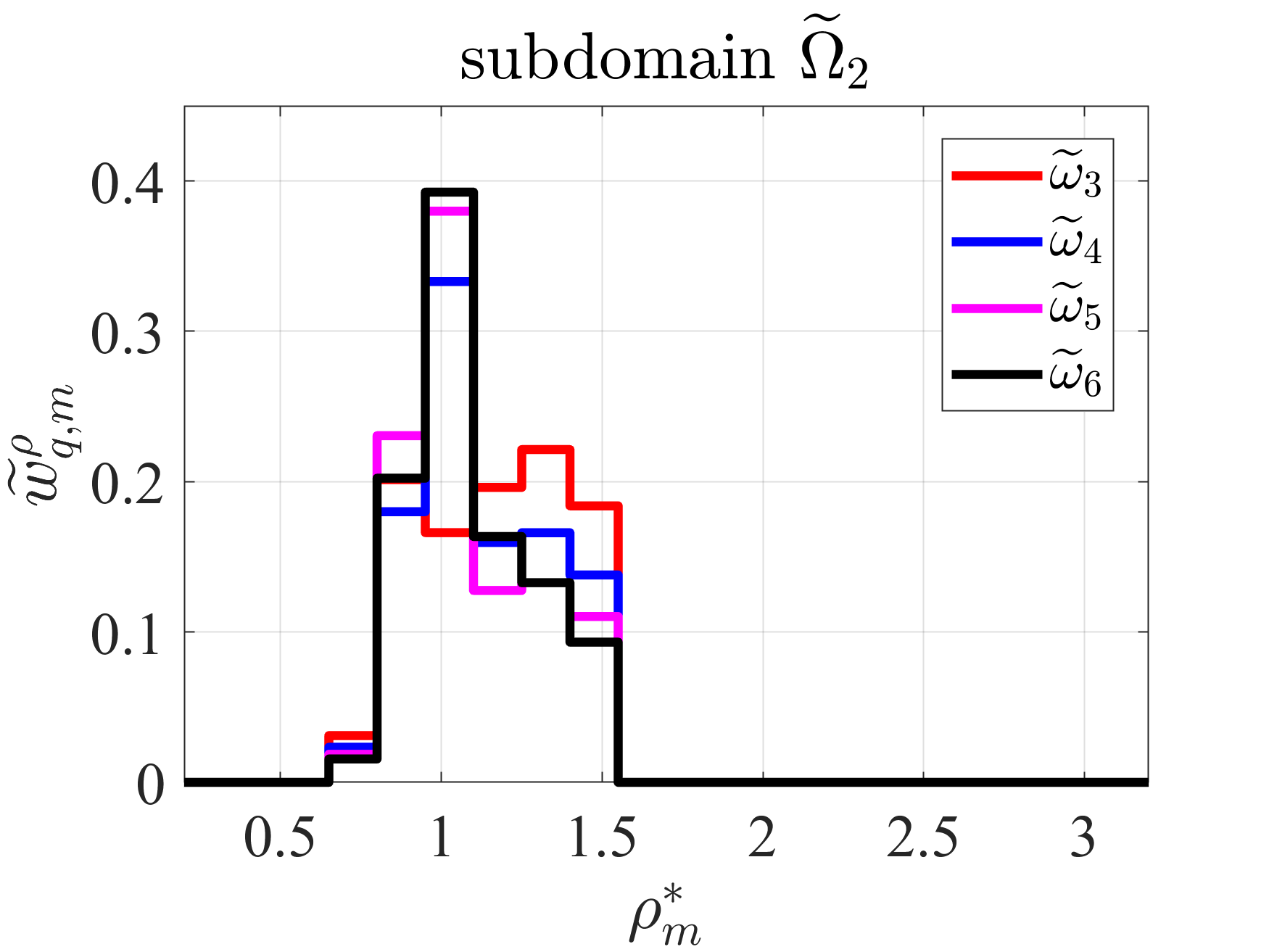}}
\caption{\sf KH Instability: Probability masses  $\widetilde w_{q,m}^{\rho}$ of the averaged density marginals $\widetilde\omega_q$, $q=3,\ldots,6$, in the subdomains
$\widetilde\Omega_1$ (left) and $\widetilde\Omega_2$ (right).\label{fig66}}
\end{figure}

The weights $\widetilde w_{q,m}^{\rho}$ are shown in Figure~\ref{fig66}. For each subdomain, the four curves represent $\widetilde\omega_3$, $\widetilde\omega_4$, $\widetilde\omega_5$, and $\widetilde\omega_6$. The height of each horizontal segment is $\widetilde w_{q,m}^{\rho}$, namely, the mass assigned to the density cell centered at $\rho_m^*$. Most of the probability mass is concentrated on several neighboring density cells, while the masses outside this region are negligible. The differences among the four averages are mainly caused by a distribution of weights among these cells. These differences indicate that the averaged marginals stabilize as higher orders are included.

\begin{figure}[ht!]
\centerline{\includegraphics[trim=0.0cm 0.cm  1.4cm 0.1cm, clip, width=7.5cm]{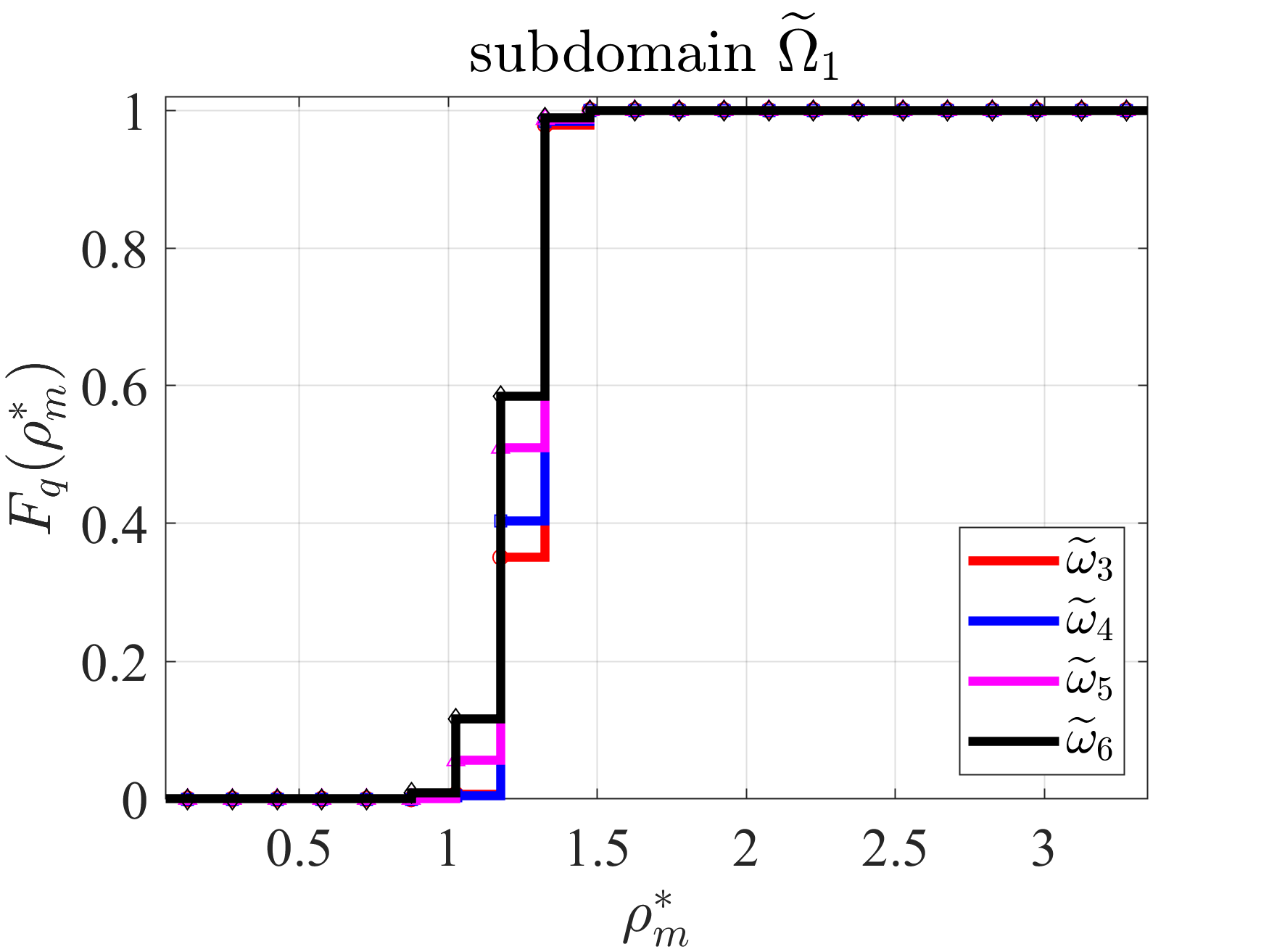}\hspace*{1.cm}
            \includegraphics[trim=0.0cm 0.0cm 1.4cm 0.1cm, clip, width=7.5cm]{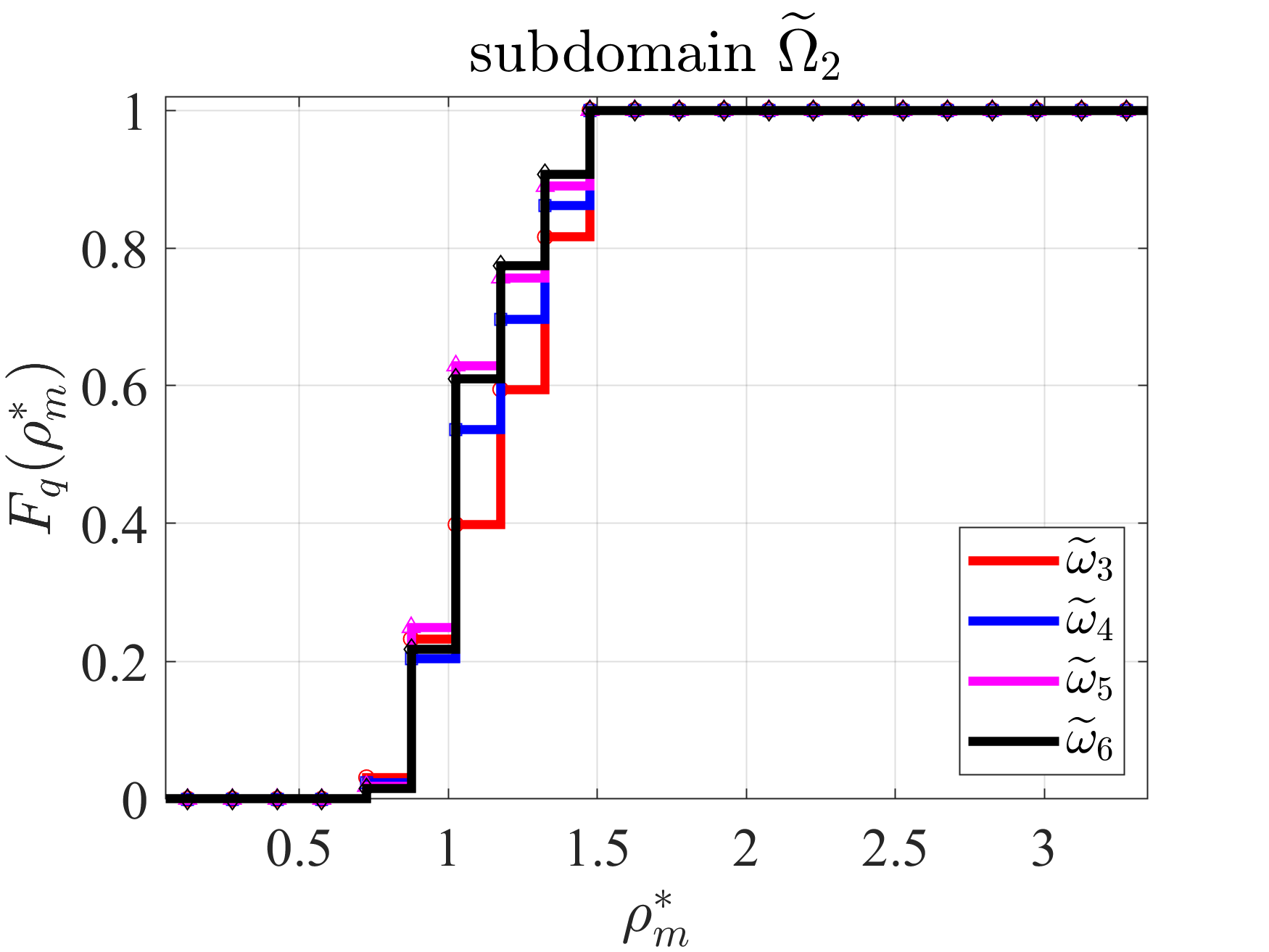}}
\caption{\sf KH Instability: Cumulative distribution functions $F_q(\rho_m^*)$ of the averaged density marginals $\widetilde\omega_q$, $q=3,\ldots,6$, in the subdomains
$\widetilde\Omega_1$ (left) and $\widetilde\Omega_2$ (right).\label{fig77}}
\end{figure}

We have also plotted the corresponding cumulative distribution functions
\begin{equation*}
 F_{q}(\rho_m^*):=\sum_{s=1}^{m}\widetilde w_{q,s}^{\rho},\quad q=3,\ldots,6,
\end{equation*}
in Figure~\ref{fig77}. Here, $F_q(\rho_m^*)$ is the probability assigned to density values not exceeding $\rho_m^*$, and its jump at $\rho_m^*$ is exactly $\widetilde w_{q,m}^{\rho}$. All four cumulative distribution functions are nondecreasing and reach one at the largest density support value, confirming the normalization of the averaged density marginals. Their closer agreement for the later averages is consistent with Figure~\ref{fig66}.

\section{Conclusions}\label{sec5}

We have extended the linear-programming-based Young-measure schemes to two space dimensions and applied them to the Kelvin--Helmholtz instability governed by the compressible Euler equations. Comparisons with the standard LLF and LLF--ELCD schemes reveal scheme-dependent differences in the small-scale vortical structures and PDFs. We have also examined several possible selection criteria for DW solutions. At every spatial reconstruction order considered, the YM schemes produce the largest time-averaged domain-integrated physical entropy. By contrast, the LLF--ELCD schemes do not systematically produce larger time-averaged domain-integrated physical entropy than the corresponding standard LLF schemes. Therefore, using entropy only to select an interface state for the LCD does not enforce a maximum-entropy selection principle. The YM results are consistent with the local maximum-entropy bias imposed through the objective function of the linear-programming problem. Moreover, the maximum-entropy, minimum-energy-defect, and minimum-Bregman-functional criteria favor the YM schemes. The finite cross-order and time averages become more stable as additional reconstruction orders are included, while the averaged density marginals remain concentrated on several neighboring density states. These fixed-mesh observations provide evidence of scheme-dependent numerical selection. Overall, the results indicate that the objective function in the linear-programming problem is an active component of the Young-measure closure rather than merely an auxiliary computational device. Future work will investigate alternative objective functions and the effects of spatial and phase-space refinement.

\begin{DA}
\paragraph{Funding.} The work of S. Chu and M. Herty was funded by the Deutsche Forschungsgemeinschaft (DFG, German Research Foundation) -- SPP 2410 Hyperbolic Balance Laws in Fluid Mechanics: Complexity, Scales, Randomness (CoScaRa) within the Project HE5386/27-2 (Zuf\"allige kompressible Euler Gleichungen: Numerik und ihre Analysis, 525853336).

\paragraph{Acknowledgments.} The authors would like to thank Prof. M\'aria Luk\'a\v{c}ov\'a-Medvi\v{d}ov\'a for insightful discussions on this work.

\paragraph{Conflicts of interest.} On behalf of all authors, the corresponding author states that there is no conflict of interest.

\paragraph{Data and software availability.} The data that support the findings of this study and FORTRAN codes developed by the authors and
used to obtain all of the presented numerical results are available from the corresponding author upon reasonable request.
\end{DA}

\appendix

\section{A Brief Review of the Reconstruction Operators}\label{appa}
In this appendix, we review the $x$-directional reconstruction operators ${\cal R}^{(r),\pm}_{\jph,k}$ at the interface $(x_{\jph},y_k)$. The $y$-directional reconstruction operators ${\cal R}^{(r),\pm}_{j,\kph}$ are constructed analogously, and their details are omitted for the sake of brevity. For further details, we refer the reader to \cite{Jiang96}.

\subsection{Local Characteristic Reconstruction}
We first define $\widehat A_{\jph,k}:=A\big(\widehat\mU_{\jph,k}\big)$, where $A(\mU):$ $=\frac{\partial\mF}{\partial\mU}$ is the $x$-directional flux Jacobian, and compute the matrices $R_{\jph,k}$ and $R^{-1}_{\jph,k}$ such that $R^{-1}_{\jph,k}$ $\widehat A_{\jph,k} R_{\jph,k}$ is diagonal. Here, $\widehat\mU_{\jph,k}$ is an average of $\mU_{j,k}$ and $\mU_{j+1,k}$, which is obtained by averaging the primitive variables:
$$
\widehat\rho_{\jph,k}=\frac{\rho_{j,k}+\rho_{j+1,k}}{2},\quad
\widehat u_{\jph,k}=\frac{u_{j,k}+u_{j+1,k}}{2},\quad
\widehat v_{\jph,k}=\frac{v_{j,k}+v_{j+1,k}}{2},\quad
\widehat p_{\jph,k}=\frac{p_{j,k}+p_{j+1,k}}{2}.
$$
The split fluxes on the reconstruction stencil are projected onto this single characteristic basis according to
\begin{equation*}
\bm Q_{\ell,k}^{x,\pm}=R^{-1}_{\jph,k}\bm F_{\ell,k}^{\pm},\qquad \ell\in{\cal S}^{(r)}_{\jph}.
\end{equation*}
The scalar left- and right-biased reconstruction operators described below are then applied to each characteristic component to obtain
\begin{equation*}
\widehat{\bm Q}_{\jph,k}^{x,+}={\cal R}^{(r),+}_{\jph}\bigg(\left\{\bm Q_{\ell,k}^{x,+}\right\}_{\ell\in{\cal S}^{(r)}_{\jph}}\bigg),\qquad
\widehat{\bm Q}_{\jph,k}^{x,-}={\cal R}^{(r),-}_{\jph}\bigg(\left\{\bm Q_{\ell,k}^{x,-}\right\}_{\ell\in{\cal S}^{(r)}_{\jph}}\bigg).
\end{equation*}
The $x$-directional numerical flux is then recovered from
\begin{equation*}
\bm{{\cal F}}^{(r)}_{\jph,k}=R_{\jph,k}\bigg(\widehat{\bm Q}_{\jph,k}^{x,+}+\widehat{\bm Q}_{\jph,k}^{x,-}\bigg).
\end{equation*}

The $y$-directional numerical fluxes $\bm{{\cal G}}^{(r)}_{j,\kph}$ are obtained analogously using the Jacobian $B(\mU):=\frac{\partial\mG}{\partial\mU}$ and its corresponding left and right eigenvector matrices. In the following subsections, the scalar formulas are applied separately to every characteristic component, and the $k$-index is suppressed for simplicity.

\begin{rmk}
For detailed explanations of how the matrices $R_{\jph,k}$ and $R^{-1}_{\jph,k}$ are computed for the 2-D Euler equations of gas dynamics, we refer the reader to \cite[Appendix~C]{CCHKL_22}.
\end{rmk}

\subsection{First-Order Reconstruction}

For $r=1$, the reconstruction operators are
\begin{equation*}
{\cal R}^{(1),+}_{\jph}\left(\{q_\ell^+\}\right)=q_j^+,
\qquad
{\cal R}^{(1),-}_{\jph}\left(\{q_\ell^-\}\right)=q_{j+1}^-.
\end{equation*}
Their sum gives a first-order LLF-type numerical flux. For the classical schemes, this is the standard LLF flux, whereas for the Young-measure scheme, the physical flux values are replaced by their Young-measure averages.

\subsection{Second-Order Reconstruction}
For $r=2$, the one-sided point values are reconstructed as
\begin{equation*}
{\cal R}^{(2),+}_{\jph}\left(\{q_\ell^+\}\right)=q_j^++\frac{\dx}{2}(q_j^+)_x,\qquad
{\cal R}^{(2),-}_{\jph}\left(\{q_\ell^-\}\right)=q_{j+1}^--\frac{\dx}{2}(q_{j+1}^-)_x,
\end{equation*}
where the slopes $(q_j^\pm)_x$ are computed using a nonlinear limiter to ensure the non-oscillatory property of the reconstruction. In this paper, we use the generalized minmod limiter (see, e.g., \cite{lie03,Nessyahu90,Sweby84}):
\begin{equation}\label{equa.6}
(q_j^\pm)_x=\operatorname{minmod}\left(\theta\frac{q_j^\pm-q_{j-1}^\pm}{\dx},~\frac{q_{j+1}^\pm-q_{j-1}^\pm}{2\dx},~\theta\frac{q_{j+1}^\pm-q_j^\pm}{\dx}\right).
\end{equation}
Here, the minmod function is defined by
\begin{equation*}
\operatorname{minmod}(z_1,z_2,\ldots):=
\begin{cases}
\min_\ell\{z_\ell\},& z_\ell>0\quad\forall \ell,\\
\max_\ell\{z_\ell\},& z_\ell<0\quad\forall \ell,\\
0,&\mbox{otherwise}.
\end{cases}
\end{equation*}
The parameter $\theta$ in \eref{equa.6} controls the amount of numerical diffusion in the resulting scheme. In general, larger values of $\theta$ lead to sharper but potentially more oscillatory reconstructions. In the numerical example reported in \S\ref{sec4}, we take $\theta=2$.

\subsection{High-Order WENO Reconstructions}
The third-, fifth-, seventh-, and ninth-order reconstruction operators are described simultaneously by setting $r=2s-1$, with $s=2,3,4,5$. Following the WENO construction, the reconstructed interface values are convex combinations of candidate values:
\begin{equation}\label{equa.17}
{\cal R}^{(2s-1),+}_{\jph}\left(\{q_\ell^+\}\right)=\sum_{m=0}^{s-1}\omega_m^+p_{m,\jph}^+,\qquad
{\cal R}^{(2s-1),-}_{\jph}\left(\{q_\ell^-\}\right)=\sum_{m=0}^{s-1}\omega_m^-p_{m,\jph}^-.
\end{equation}
The nonlinear weights are given by 
$\alpha_m^\pm=\frac{d_m^{(s)}}{(\varepsilon_0+\beta_m^\pm)^2}$, where $\omega_m^\pm={\alpha_m^\pm}\bigg/{\displaystyle\sum_{\ell=0}^{s-1}\alpha_\ell^\pm}$, and $\varepsilon_0>0$ is a fixed small parameter. In this paper, we take $\varepsilon_0=10^{-12}$ in all of the schemes. The optimal weights are
\begin{equation*}
d_m^{(s)}=\frac{\binom{s-1}{m}\binom{s}{m}}{\binom{2s-1}{s}},\qquad m=0,\ldots,s-1.
\end{equation*}
For the fifth-, seventh-, and ninth-order reconstructions, the smoothness indicators are computed by
\begin{equation}\label{equa.13}
\beta_m^+=\sum_{\nu=1}^{s-1}\dx^{2\nu-1}\int_{I_j}\left[\frac{{\rm d}^\nu}{{\rm d}x^\nu}p_m^+(x)\right]^2{\rm d}x,\qquad
\beta_m^-=\sum_{\nu=1}^{s-1}\dx^{2\nu-1}\int_{I_{j+1}}\left[\frac{{\rm d}^\nu}{{\rm d}x^\nu}p_m^-(x)\right]^2{\rm d}x,
\end{equation}
where $I_j=[x_{j-1/2},x_{j+1/2}]$, and the polynomials $p_m^\pm(x)$ have degree at most $s-1$ and satisfy
\begin{equation*}
\frac{1}{\dx}\int_{x_{\ell-1/2}}^{x_{\ell+1/2}}p_m^\pm(x)\,{\rm d}x=q_\ell^\pm
\end{equation*}
for every index $\ell$ in the corresponding candidate stencil. This is the standard auxiliary cell-average interpretation used to construct finite-difference WENO numerical fluxes from the pointwise split-flux data $q_\ell^\pm$. We refer to \cite{Gao20} for the explicit forms of the smoothness indicators used in the fifth-, seventh-, and ninth-order schemes. For the third-order scheme, the indicators in \eref{equa.13} are replaced by the new smoothness indicators introduced in \cite[Appendix ~A]{CHT25} to achieve higher resolution, while the remaining reconstruction formulas are unchanged.

Finally, the interface values in \eref{equa.17} are given by
\begin{equation*}
p_{m,\jph}^+=\sum_{\ell=0}^{s-1}C_{m,\ell}^{[s]}q_{j-s+1+m+\ell}^+,\qquad
p_{m,\jph}^-=\sum_{\ell=0}^{s-1}C_{m,\ell}^{[s]}q_{j+s-m-\ell}^-,
\end{equation*}
where the coefficient matrices are given by
\begin{equation*}
\bm C^{[2]}=\begin{pmatrix}
-\dfrac12&\dfrac32\\[3mm]
\dfrac12&\dfrac12
\end{pmatrix},\qquad
\bm C^{[3]}=
\begin{pmatrix}
\dfrac13&-\dfrac76&\dfrac{11}{6}\\[3mm]
-\dfrac16&\dfrac56&\dfrac13\\[3mm]
\dfrac13&\dfrac56&-\dfrac16
\end{pmatrix},\qquad
\bm C^{[4]}=
\begin{pmatrix}
-\dfrac14&\dfrac{13}{12}&-\dfrac{23}{12}&\dfrac{25}{12}\\[3mm]
\dfrac1{12}&-\dfrac5{12}&\dfrac{13}{12}&\dfrac14\\[3mm]
-\dfrac1{12}&\dfrac7{12}&\dfrac7{12}&-\dfrac1{12}\\[3mm]
\dfrac14&\dfrac{13}{12}&-\dfrac5{12}&\dfrac1{12}
\end{pmatrix},
\end{equation*}
and
\begin{equation*}
\bm C^{[5]}=
\begin{pmatrix}
\dfrac15&-\dfrac{21}{20}&\dfrac{137}{60}
&-\dfrac{163}{60}&\dfrac{137}{60}\\[3mm]
-\dfrac1{20}&\dfrac{17}{60}&-\dfrac{43}{60}
&\dfrac{77}{60}&\dfrac15\\[3mm]
\dfrac1{30}&-\dfrac{13}{60}&\dfrac{47}{60}
&\dfrac9{20}&-\dfrac1{20}\\[3mm]
-\dfrac1{20}&\dfrac9{20}&\dfrac{47}{60}
&-\dfrac{13}{60}&\dfrac1{30}\\[3mm]
\dfrac15&\dfrac{77}{60}&-\dfrac{43}{60}
&\dfrac{17}{60}&-\dfrac1{20}
\end{pmatrix}.
\end{equation*}

\section{Entropy-Based Local Characteristic Decomposition}\label{appb}
In this appendix, we briefly describe the entropy-based LCD; see, e.g., \cite{CHK_NEW_LCD}. The reconstruction procedure is the same as that presented in Appendix~\ref{appa}, except that the interface state $\widehat\mU_{\jph,k}$ used to construct the matrix $\widehat A_{\jph,k}=A(\widehat\mU_{\jph,k})$ is selected according to the physical entropy,
determined by the following algorithm.
\begin{algorithm}[H]
\caption{Entropy-based average interface state in the $x$-direction.}
\begin{algorithmic}[1]
\State {\bf Set} ${S}_1:={S}(\,\rho_{j,k},p_{j,k})$, ${S}_2:={S}(\,\rho_{j+1,k},p_{j+1,k})$,
${S}_3:={S}(\,\rho_{j,k},p_{j+1,k})$, ${S}_4:={S}(\,\rho_{j+1,k}, p_{j,k})$
\vskip3pt
\If{${S}_1=\max\{{S}_1, {S}_2, {S}_3, {S}_4\}$}
\State set $\widehat\rho_{\jph,k}:=\rho_{j,k}$,
$\widehat u_{\jph,k}:=u_{j,k}$,
$\widehat v_{\jph,k}:=v_{j,k}$, 
$\widehat p_{\jph,k}:=p_{j,k}$
\ElsIf{${S}_2=\max\{{S}_1, {S}_2, {S}_3, {S}_4\}$}
\State set $\widehat\rho_{\jph,k}:=\rho_{j+1,k}$,
$\widehat u_{\jph,k}:=u_{j+1,k}$,
$\widehat v_{\jph,k}:=v_{j+1,k}$, 
$\widehat p_{\jph,k}:=p_{j+1,k}$
\ElsIf{${S}_3=\max\{{S}_1, {S}_2, {S}_3, {S}_4\}$}
\State set $\widehat\rho_{\jph,k}:=\rho_{j,k}$,
$\widehat u_{\jph,k}:=(u_{j,k}+u_{j+1,k})/2$,
$\widehat v_{\jph,k}:=(v_{j,k}+v_{j+1,k})/2$,
$\widehat p_{\jph,k}:=p_{j+1,k}$
\Else
\State set $\widehat\rho_{\jph,k}:=\rho_{j+1,k}$,
$\widehat u_{\jph,k}:=(u_{j,k}+u_{j+1,k})/2$,
$\widehat v_{\jph,k}:=(v_{j,k}+v_{j+1,k})/2$, 
$\widehat p_{\jph,k}:=p_{j,k}$
\EndIf
\end{algorithmic}
\end{algorithm}
The remaining reconstruction steps are exactly the same as those described in Appendix~\ref{appa}. The entropy-based interface state in the $y$-direction is constructed analogously, and its details are omitted for the sake of brevity.

\medskip 
\bibliographystyle{siamnodash}
\bibliography{reference}
\end{document}